\documentclass[12pt,twoside]{article}
\usepackage[english]{babel}
\usepackage{umj_3_1}           % Style package of Ural Mathematical Journal
\usepackage{graphicx, srcltx}   % For graphics
\usepackage{cite}               % For correct bibliography
\usepackage{url}                                % use command \cite{1,2,3,4,8,9,10}
\usepackage{amsmath}            % If the package features are used
\usepackage{amsfonts,amssymb}   % If there is need to use math alphabets
\usepackage{srcltx}             % For easy move between WiNEdT and YAP

\newcommand{\Leftclt}{Victor Nijimbere}   % Left page header - name & family name of article author
\newcommand{\Rightclt}{Some non-elementary integrals of sine, cosine and exponential
integrals type} %% Right page header - Short title of article
\usepackage{authblk}

\begin{document}

\date{}

\Mainclt            % Page header print on the first page

% Enter article title and authors
% Reference to grants is a footnote to article title. You may use also Acknowledgments section for this purpose if you wish.

\begin{Titul}
%{\large \bf Evaluation of some non-elementary integrals involving sine, cosine, exponential and logarithmic
%integrals: Part I}\\ \vspace{5mm}
{\large \bf EVALUATION OF SOME NON-ELEMENTARY INTEGRALS INVOLVING SINE, COSINE, EXPONENTIAL AND LOGARITHMIC INTEGRALS: PART I}\\ \vspace{5mm}
%\footnote{This work was supported by ... (project no.NN-N). 
%You may use also Acknowledgments section for this purpose if you wish.}\\[3ex]
{{\bf Victor Nijimbere} \\ [2ex] {\small School of Mathematics and Statistics, Carleton University, \hspace{7cm}Ottawa, Ontario, Canada,\hspace{9cm}victornijimbere@gmail.com}\\[3ex]} 
%{\small %victornijimbere@gmail.com}\\ [4ex]}
%{\bf  Second B. Author}\\[2ex] {\small University, City, Country, Email }\\[5ex]}
\end{Titul}

\begin{Anot}
{\bf Abstract:} The non-elementary integrals $\text{Si}_{\beta,\alpha}=\int [\sin{(\lambda x^\beta)}/(\lambda x^\alpha)] dx,\beta\ge1,\alpha\le\beta+1$ and $\text{Ci}_{\beta,\alpha}=\int [\cos{(\lambda x^\beta)}/(\lambda x^\alpha)] dx, \beta\ge1, \alpha\le2\beta+1$, where $\{\beta,\alpha\}\in\mathbb{R}$, are evaluated in terms of the hypergeometric functions $_{1}F_2$ and $_{2}F_3$, and their asymptotic expressions for $|x|\gg1$ are also derived. The integrals of the form $\int [\sin^n{(\lambda x^\beta)}/(\lambda x^\alpha)] dx$ and $\int [\cos^n{(\lambda x^\beta)}/(\lambda x^\alpha)] dx$, where $n$ is a positive integer, are expressed in terms $\text{Si}_{\beta,\alpha}$ and $\text{Ci}_{\beta,\alpha}$, and then evaluated.  $\text{Si}_{\beta,\alpha}$ and $\text{Ci}_{\beta,\alpha}$ are also evaluated in terms of the hypergeometric function $_{2}F_2$. And so, the hypergeometric functions, $_{1}F_2$ and $_{2}F_3$, are expressed in terms of $_{2}F_2$.
The exponential integral $\text{Ei}_{\beta,\alpha}=\int (e^{\lambda x^\beta}/x^\alpha) dx$ where $\beta\ge1$ and $\alpha\le\beta+1$ and the logarithmic integral $\text{Li}=\int_{\mu}^{x} dt/\ln{t}, \mu>1$ are also expressed in terms of $_{2}F_2$, and their asymptotic expressions are investigated. It is found that for $x\gg\mu$, $\text{Li}\sim {x}/{\ln{x}}+\ln{\left(\frac{\ln{x}}{\ln{\mu}}\right)}-2-
\ln{\mu}\hspace{.075cm} _{2}F_{2}(1,1;2,2;\ln{\mu})$, where the term $\ln{\left(\frac{\ln{x}}{\ln{\mu}}\right)}-2-
\ln{\mu}\hspace{.075cm} _{2}F_{2}(1,1;2,2;\ln{\mu})$ is added to the known expression in mathematical literature $\text{Li}\sim {x}/{\ln{x}}$.
The method used in this paper consists of expanding the integrand as a Taylor and integrating the series term by term, and can be used to evaluate the other cases which are not considered here. This work is motivated by the applications of sine, cosine exponential and logarithmic integrals in Science and Engineering, and some applications are given.

\end{Anot}

%%%%%%%%%%%%%%%%%%%%%%%%%%%%%%%%%%%%%%%%%%%%
 %Please replace what follows by the body %%
 %of your article                         %%
                                          %%
%%%%%%%%%%%%%%%%%%%%%%%%%%%%%%%%%%%%%%%%%%%%

\section{Introduction}\label{sec:1}
\setcounter{equation}{0}

%This article illustrates the preparation of a paper for Ural Mathematical Journal using \LaTeX2e.

%Your paper must be prepared in {\bf English}. To ensure that your paper will be reproduced
%clearly and in the proper size and form, please observe the following instructions.

%The length of the title of your paper should not exceed two lines. 
%The paper should begin with an abstract, followed by key words, and then an introduction.
%At the end, it should feature a conclusion and a list of references.
%The maximum length of your paper, including figures and tables,should be 10-20 pages.

%In addition to length and
%format restrictions given here, all papers will be reviewed and must meet technical standards.

%\begin{verbatim}\begin{assen}...\end{assen}\end{verbatim}
%For example,
%\begin{teoen}[{\cite[Theorem~2]{DeVore-Lorentz, Geronimus-1936}}]
%Theorem text goes here.
%\end{teoen}
%\begin{teoen}[(Eremenko, Yuditskii {\rm \cite[p.~25]{Eremenko-Yuditskii}})]
%Theorem text goes here.
%\end{teoen}
%\proofen    Proof text goes here.  \hfill$\square$\\[1ex]%--- P r o o f.
%\proofnpen follows from $\ldots$\,... \hfill$\square$\\[1ex]% --- P r o o f

%\emph{Remark~1.} Text...\\[1ex] % Remark (Çàìå÷àíèå)
%\emph{Example~1.} Text...\\[1ex] % Example

%We use \hfill$\square$ for the end of the proof.

\begin{defen} An elementary function is a function of one variable constructed using that variable and constants, and by performing a finite number of repeated algebraic operations involving exponentials and logarithms.
An indefinite integral which can be expressed in terms of elementary functions is an elementary integral. And if, on the other hand, it cannot be evaluated in terms of elementary functions, then it is non-elementary \cite{MZ,R}.
\label{defen1}
\end{defen}

Liouville 1938's Theorem gives conditions to determine whether  a given integral is elementary or non-elementary \cite{MZ,R}. For instance, it was shown in \cite {MZ,R}, using  Liouville 1938's Theorem, that the integral $\text{Si}_{1,1}=\int (\sin{x}/x) dx$ is non-elementary. With similar arguments as in \cite {MZ,R}, One can show that $\text{Ci}_{1,1}=\int (\cos{x}/x) dx$ is also non-elementary. Using the Euler formulas $e^{\pm ix}=\cos{x}\pm i\sin{x}$, and noticing that if the integral of a function $g(x)$ is elementary, then both its real and imaginary parts are elementary \cite{MZ}, one can, for instance, prove that the integrals $\text{Si}_{\beta,\alpha}=\int [\sin{(\lambda x^\beta)}/(\lambda x^\alpha)] dx, \beta\ge1, \alpha\ge1$, and $\text{Ci}_{\beta,\alpha}=\int [\cos{(\lambda x^\beta)}/(\lambda x^\alpha)] dx$, where $\beta\ge1$ and $\alpha\ge1$, are non-elementary by using the fact that their real and imaginary parts are non-elementary.
The integrals $\int [\sin^n{(\lambda x^\beta)}/(\lambda x^\alpha)] dx$ and $\int [\cos^n{(\lambda x^\beta)}/(\lambda x^\alpha)] dx$, where $n$ is a positive integer, are also non-elementary since they can be expressed in terms of $\text{Si}_{\beta,\alpha}$ and $\text{Ci}_{\beta,\alpha}$.

To my knowledge, no one has evaluated these integrals before. To this end, in this paper, formulas for these non-elementary integrals are expressed in terms of the hypergeometric functions  $_{1}F_2$ and $_{2}F_3$ 
whose properties, for example, the asymptotic expansions for large argument ($|\lambda x|\gg1$), are known \cite{N}. We do so by expanding the integrand in terms of its Taylor series and by integrating the series term by term as in \cite{NV}. And therefore, their corresponding definite integrals can be evaluated using the Fundamental Theorem of Calculus (FTC). For example, the sine integral
\begin{equation}
\text{Si}_{\beta,\alpha}=\int\limits_{A}^{B} \frac{\sin {(\lambda x^\beta)}}{(\lambda x^\alpha)} dx,\hspace{.2cm} \beta\ge1,,\hspace{.2cm} \alpha\le\beta+1,
\label{eq1}
\end{equation}
is evaluated for any $A$ and $B$ using the FTC.

On the other hand, the integrals  $\text{Ei}_{\beta,\alpha}=\int (e^{\lambda x^\beta}/x^\alpha) dx$ and $\int dx/\ln{x}$, are expressed in terms of the hypergeometric function $_{2}F_2$. This is quite important since one may re-investigate the asymptotic behavior of the exponential ($\text{Ei}$) and logarithmic ($\text{Li}$) integrals \cite{H} using the asymptotic expressions of the hypergeometric function $_{2}F_2$ which are known \cite{N}.

Some other non-elementary integrals which can be written in terms of $\text{Ei}_{\beta,\alpha}$ or $\int dx/\ln{x}$ are also evaluated.  For instance, as a result of substitution, the integral $\int e^{\lambda e^{\beta x}} dx$ is written in terms of $\text{Ei}_{\beta,1}=\int (e^{\lambda x^\beta}/x) dx$ and then evaluated in terms of $_{2}F_2$, and using integration by parts, the integral $\int \ln(\ln{x}) dx$ is written in terms of $\int dx/\ln{x}$ and then evaluated in terms of $_{2}F_2$ as well.

Using the Euler identity $e^{\pm ix}=\cos(x)\pm i\sin(x)$ or the hyperbolic identity $e^{\pm x}=\cosh(x)\pm\sinh(x)$, $\text{Si}_{\beta,\alpha}$ and $\text{Ci}_{\beta,\alpha}$ are evaluated in terms $\text{Ei}_{\beta,\alpha}$. And hence, the hypergeometric functions $_{1}F_2$ and $_{2}F_3$ are expressed in terms of the hypergeometric $_{2}F_2$.

This type of integrals find applications in many fields in Science and Engineering. For instance, in wireless telecommunications,  the random attenuation capacity of a channel, known as fading capacity, is calculated as \cite{MM} 
\begin{equation}
C_{\mbox{fading}}=E[\log_2(1+P|H|^2)]=\int\limits_0^\infty \log_2(1+P\xi)e^{-\xi} d\xi=\frac{1}{\ln{2}}e^{\frac{1}{p}}
\left[E_{1,1}\left(\infty\right)-E_{1,1}\left(\frac{1}{p}\right)\right],
\label{application1}
\end{equation}
where the fading coefficient $H$ is a complex Gaussian random variable, and  $E\left(|X|^2\le P\right)$ is the maximum average transmitted power of a complex-valued channel input $X$.
In number theory, the prime number Theorem states that \cite{H}
\begin{equation}
\pi(x)\sim\mbox{Li}(x)=\int\limits_\mu^x \frac{dx}{\ln{x}}, \mu>1,
\label{application2}
\end{equation}
where $\pi(x)$ denotes the number of primes small than or equal to $x$. Moreover, there are applications of sine and cosine integrals in  electromagnetic theory, see for example Lebedev \cite{L}. Therefore, it is quite important to adequately evaluate these integrals. 

For that reason, the main goal of this paper is to evaluate non-elementary integrals of sine, cosine, exponential and logarithmic integrals type in terms of elementary and special functions with well known properties so that the fundamental theorem of calculus can be used so that we can avoid to use numerical integration. 

Part I is indeed devoted to the cases $\text{Si}_{\beta,\alpha}=\int [\sin{(\lambda x^\beta)}/(\lambda x^\alpha)] dx,\beta\ge1,\alpha\le\beta+1$, $\text{Ci}_{\beta,\alpha}=\int [\cos{(\lambda x^\beta)}/(\lambda x^\alpha)] dx, \beta\ge1, \alpha \le2\beta+1$ and  $\text{Ei}_{\beta,\alpha}=\int (e^{\lambda x^\beta}/x^\alpha) dx$ where $\beta\ge1, \alpha\le\beta+1$, where $\{\beta,\alpha\}\in\mathbb{R}$.  The other cases  $\text{Si}_{\beta,\alpha}=\int [\sin{(\lambda x^\beta)}/(\lambda x^\alpha)] dx,\beta\ge1,\alpha>\beta+1$, $\text{Ci}_{\beta,\alpha}=\int [\cos{(\lambda x^\beta)}/(\lambda x^\alpha)] dx, \beta\ge1, \alpha>2\beta+1$ and  $\text{Ei}_{\beta,\alpha}=\int (e^{\lambda x^\beta}/x^\alpha) dx$ where $\beta\ge1,\alpha>\beta+1$, where $\{\beta,\alpha\}\in\mathbb{R}$, which may involve series whose properties are not necessary known will be considered in Part 2 \cite{NV2}.

 Before we proceed to the objectives of this paper (see sections \ref{sec:2},  \ref{sec:3},  \ref{sec:4} and \ref{sec:5}), we first define the generalized hypogeometric function as it is an important mathematical that we are going to use throughout the paper.

\begin{defen}The generalized hypergeometric function, denoted as $_pF_q$, is a special function given by the series \cite{AS,N}
\begin{equation}
_p F_q(a_1, a_2,\cdots,a_p;b_1, b_2, \cdots, b_q; x)=\sum\limits_{n=0}^{\infty}\frac{(a_1)_n (a_2)_n\cdots (a_p)_n}{(b_1)_n (b_2)_n\cdots (b_q)_n}\frac{x^n}{n!},
\label{hypergeometric}
\end{equation} 
where $a_1, a_2,\cdots,a_p$ and $;b_1, b_2, \cdots, b_q$ are arbitrary constants, $(\vartheta)_n=\Gamma(\vartheta+n)/\Gamma(\vartheta)$ (Pochhammer's notation \cite{1}) for any complex $\vartheta$, with $(\vartheta)_0=1$, and $\Gamma$ is the standard gamma function \cite{AS,N}.
\label{defen2}
\end{defen}

\section{Evaluation of the sine integral and related integrals}\label{sec:2}
\setcounter{equation}{0}
%Symbols and acronyms should be defined the first time they appear. Use the International System (SI) of units.

\begin{prpen}
The function $G(x)={x}\hspace{.075cm} _{1}F_{2}\left(\frac{1}{2};\frac{3}{2},\frac{3}{2};-\frac{\lambda^2 x^2}{4}\right)$, where ${}_1F_2$ is a hypergeometric function \cite{AS} and $\lambda$ is an arbitrarily constant, is the antiderivative of the function $g(x)=\frac{\sin{(\lambda x)
}}{\lambda x}$. Thus,
\begin{equation}
\int\frac{\sin{(\lambda x)}}{\lambda x}dx={x}\hspace{.075cm} _{1}F_{2}\left(\frac{1}{2};\frac{3}{2},\frac{3}{2};-\frac{\lambda^2 x^2}{4}\right)+C.
\label{eq2}
\end{equation}
\label{prpen:1}
\end{prpen}

\proofen   To prove Proposition \ref{prpen:1},
we expand $g(x)$ as Taylor series and integrate the series term by term.
We also use the gamma duplication formula \cite{AS}
\begin{equation}
\Gamma(2\alpha)=(2\pi)^{-\frac{1}{2}}2^{2\alpha-\frac{1}{2}}\Gamma(\alpha)\Gamma\left(\alpha+\frac{1}{2}\right), \hspace{.15cm} \alpha\in\mathbb{C},
\label{eq3}
\end{equation}
the Pochhammer's notation for the gamma function \cite{AS},
\begin{equation}
(\alpha)_n=\alpha (\alpha+1)\cdot\cdot\cdot(\alpha+n-1)=\frac{\Gamma(\alpha+n)}{\Gamma(\alpha)}, \hspace{.15cm} \alpha\in\mathbb{C},
\label{eq4}
\end{equation}
and the property of the gamma function $\Gamma(\alpha+1)=\alpha\Gamma(\alpha)$ (eg., $\Gamma\left(n+\frac{3}{2}\right)=\left(n+\frac{1}{2}\right)\Gamma\left(n+\frac{1}{2}\right)$ for any real $n$). We then obtain

\begin{align}
\int g(x)dx&=\int \frac{\sin{(\lambda x)}}{\lambda x}dx
=\int\frac{1}{\lambda x}\sum\limits_{n=0}^{\infty}(-1)^n\frac{(\lambda x)^{2n+1}}{(2n+1)!}dx
%\nonumber\\ &=\lambda\int\sum\limits_{n=0}^{\infty}(-1)^n\frac{(\lambda x)^{2n}}{(2n+1)!}dx
%\nonumber\\&=\lambda\sum\limits_{n=0}^{\infty}(-1)^n\frac{\lambda^{2n}}{(2n+1)!}\int x^{2n} dx
=\lambda\sum\limits_{n=0}^{\infty}(-1)^n\frac{\lambda^{2n}}{(2n+1)!}\frac{x^{2n+1}}{2n+1}+C
\nonumber\\&=\frac{x}{2}\sum\limits_{n=0}^{\infty}(-1)^n\frac{\lambda^{2n}}{(2n+1)!}\frac{x^{2n}}{n+\frac{1}{2}}+C
=\frac{x}{2}\sum\limits_{n=0}^{\infty}\frac{\Gamma\left(n+\frac{1}{2}\right)}{\Gamma(2n+2)\Gamma\left(n+\frac{3}{2}\right)}(-\lambda^2 x^{2})^{n}+C
\nonumber\\&={x}\sum\limits_{n=0}^{\infty}\frac{\left(\frac{1}{2}\right)_n}{\left(\frac{3}{2}\right)_n\left(\frac{3}{2}\right)_n}\frac{\left(-\frac{\lambda^2 x^2}{4}\right)^{n}}{n!}+C
={x}\hspace{.075cm} _{1}F_{2}\left(\frac{1}{2};\frac{3}{2},\frac{3}{2};-\frac{\lambda^2 x^2}{4}\right)+C =
G(x)+C.
\label{eq5}
\end{align}
 \hfill$\square$\\[1ex]%--- P r o o f.

In the following Lemma, we assume that the function $G(x)$ is unknown and therefore we establish its properties such as the inflection points and its  behaviour as $x\to\pm\infty$. 

\begin{lemen}
Let $G(x)$ be the antiderivative for $g(x)=\frac{\sin{x}}{x}$ ($\lambda=1$), and $G(0)=0$.
\begin{enumerate}
\item Then $G(x)$ is linear around $x=0$ and the point $(0,G(0))=(0,0)$ is an inflection point of the curve $Y=G(x)$, $x\in\mathbb{R}$.
\item And $\lim\limits_{x\rightarrow-\infty}G(x)=-\theta$ while $\lim\limits_{x\rightarrow+\infty}G(x)=\theta$, where $\theta$ is a positive finite constant.
\end{enumerate}
\label{lemen:1}
\end{lemen}

\proofen 
\begin{enumerate}
\item %$\lim\limits_{x\rightarrow0}g(x)=1$.
 The series $g(x)=\frac{\sin{x}}{x}=\sum\limits_{n=0}^{\infty}(-1)^n\frac{(\lambda x)^{2n}}{(2n+1)!}$ gives $G^\prime(0)=g(0)=1$. Then, around $x=0$, $G(x)\sim x$ since $G^\prime(0)=g(0)=1$ and $G(0)=0$. Moreover, $G^{\prime\prime}(x)=g^\prime(x)=\left(\frac{\sin{x}}{x}\right)^\prime=-\lambda^2 x \sum\limits_{n=0}^{\infty}(-1)^n\frac{(2n+2)(\lambda x)^{2n}}{(2n+3)!}$, and so $G^{\prime\prime}(0)=g^\prime(0)=0$. Hence, by the second derivative test, the point $(0,G(0))=(0,0)$ is an inflection point of the curve $Y=G(x)$.

\item It is straight forward, using Squeeze theorem, to obtain $\lim\limits_{x\rightarrow-\infty}g(x)=\lim\limits_{x\rightarrow+\infty}g(x)=0$. And since both $g(x)$ and $G(x)$ are analytic on $\mathbb{R}$, then $G(x)$ has to be constant as $x\rightarrow\pm\infty$ by Liouville Theorem (section 3.1.3 in \cite{K}) since if a complex function is entire on $\mathbb{C}$ then both its imaginary and real parts are analytic on the real line $\mathbb{R}$ including at $x\to\pm\infty$. 
%(asymptotic expansion in Lemma \ref{lempr:1}). 
Also, there exists some numbers $\delta>0$ and $\epsilon$ such that if $|x|>\delta$ then $||\sin{x}|/x-{1}/{x}|<\epsilon$, and $\lim\limits_{x\rightarrow-\infty}(|\sin{x}|/x)/({1}/{x})=\lim\limits_{x\rightarrow+\infty}(|\sin{x}|/x)/({1}/{x})=\pm1$.
This makes the function $g_1(x)=-{1}/{x}$ an envelop of $g(x)$ away from $x=0$ if $\sin{x}<0$ and $g_2(x)={1}/{x}$ an envelop of $g(x)$ away from $x=0$ if $\sin{x}>0$. Moreover, on one hand, $g_2^\prime\le G^{\prime\prime}\le g_1^\prime$ if $x<-\delta$, and $g_1^\prime$ and $g_2^\prime$ do not change signs. While on another hand, $g_1^\prime\le G^{\prime\prime}\le g_2^\prime$ if $x>\delta$, and also $g_1^\prime$ and $g_2^\prime$ do not change signs. Therefore there exists some number $\theta>0$ such $G(x)$ oscillates about $\theta$ if $x>\delta$ and $G(x)$ oscillates about $-\theta$ if $x<-\delta$. And $|G(x)|\le\theta$ if $|x|\le\delta$.
\end{enumerate}
\hfill$\square$\\[1ex]

\begin{figure}[t]
\centerline{\includegraphics[width=0.45\textwidth]{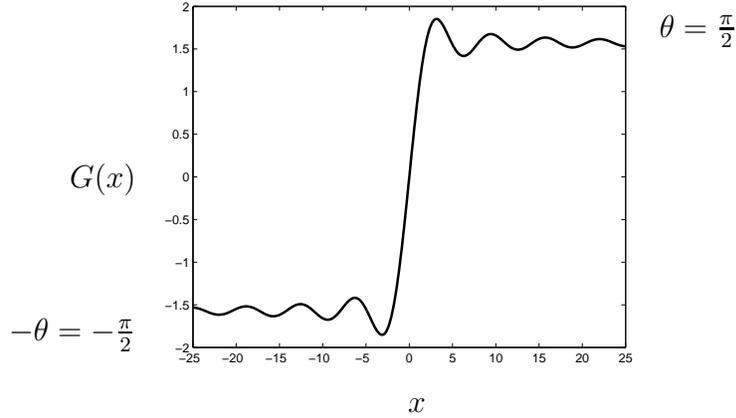}}
\begin{picture}(0,0)
\put (110,93){$G(x)$}
\put (238,8){$x$}
\put (87,35){$-\theta=-\frac{\pi}{2}$}
\put (333,150){$\theta=\frac{\pi}{2}$}
\end{picture}
\caption{$G(x)$ is the antiderivative of $\sin(x)/x$ given in (\ref{eq6}).}
\label{fig1}
\end{figure}

\emph{Example~1.}
%\begin{exaen}
For instance, if $\lambda=1$, then
\begin{equation}
\int \frac{\sin{x}}{x} dx={x}\hspace{.075cm} _{1}F_{2}\left(\frac{1}{2};\frac{3}{2},\frac{3}{2};-\frac{x^2}{4}\right)+C.
\label{eq6}
\end{equation}
%\label{ex:1}
By Proposition \ref{prpen:1}, the antiderivative of $g(x)= \frac{\sin{x}}{x}$ is $G(x)={x}\hspace{.075cm} _{1}F_{2}\left(\frac{1}{2};\frac{3}{2},\frac{3}{2};-\frac{x^2}{4}\right)$, and the graph of $G(x)$ is shown in Figure \ref{fig1}. It is in agreement with Lemma \ref{lemen:1}. It is seen in Figure \ref{fig1} that $(0,G(0))=(0,0)$ is an inflection point and that $G$ attains some constants as $x\rightarrow\pm\infty$ as predicted by Lemma \ref{lemen:1}.
%\end{exaen}
\\[1ex]

In the following lemma, we obtain the values of $G(x)$, the antiderivative of the function $g(x) =\sin{(\lambda x)}/(\lambda x)$, as $x\to\pm\infty$ using the asymptotic expansion of the hypergeometric function $_1F_2$.

\begin{lemen} Consider $G(x)$ in Proposition \ref{prpen:1},and preferably assume that $\lambda>0$.
\begin{enumerate}
\item Then,
\begin{equation}
G(-\infty)=\lim_{x\rightarrow-\infty}G(x)=\lim_{x\rightarrow\infty}x\hspace{.1cm}{}_1F_2\left(\frac{1}{2};\frac{3}{2},\frac{3}{2};-\frac{\lambda^2 x^2}{4}\right)=-\frac{{\pi}}{2\lambda},
\label{eq7}
\end{equation}
and
\begin{equation}
G(+\infty)=\lim_{x\rightarrow+\infty}G(x)=\lim_{x\rightarrow\infty}x\hspace{.1cm}{}_1F_2\left(\frac{1}{2};\frac{3}{2},\frac{3}{2};-\frac{\lambda^2 x^2}{4}\right)=\frac{{\pi}}{2\lambda}.
\label{eq8}
\end{equation}
\item And by the FTC,
\begin{equation}
\int\limits_{-\infty}^{\infty} \frac{\sin{(\lambda x)
}}{\lambda x} dx=G(+\infty)-G(-\infty)=\frac{{\pi}}{2\lambda}-\left(-\frac{\sqrt{\pi}}{2\lambda}\right)=\frac{{\pi}}{\lambda}.
\label{eq9}
\end{equation}
\end{enumerate}
\label{lemen:2}
\end{lemen}

\proofen
\begin{enumerate}
\item To prove (\ref{eq7}) and (\ref{eq8}), we use the asymptotic formula for the hypergeometric function ${}_1F_2$ which is valid for $|z|\gg 1$ and $-\pi\le\mbox{arg}\hspace{.1cm} z\le\pi$, where $\mbox{arg}\hspace{.1cm} z$ is the argument of $z$ in the complex plane. It can be derived using formulas 16.11.1, 16.11.2 and 16.11.8 in \cite{N} and is given by
\begin{align}
&{}_1F_2\left(a_1;b_1,b_2;-z\right)=\nonumber\\&\Gamma(b_1)\Gamma(b_2)z^{-a_1}\left\{\sum\limits_{n=0}^{R-1}\frac{(a_1)_n }{\Gamma(b_1-a_1-n)\Gamma(b_2-a_1-n)}\frac{(-z)^{-n}}{n!}+O(|z|^{-R})\right\}\nonumber\\&+\frac{\Gamma(b_1)\Gamma(b_2)}{\Gamma(a_1)}
+\frac{e^{2z^{\frac{1}{2}}e^{-i\frac{\pi}{2}}}(ze^{-i\pi})^{\frac{a_1-b_1-b_2+\frac{1}{2}}{2}}}{\sqrt{\pi}}\left\{\sum\limits_{n=0}^{S-1}\frac{\mu_n}{2^{n+1}}(ze^{-i\pi})^{-n}+O(|z|^{-S})\right\}
\nonumber\\&+\frac{\Gamma(b_1)\Gamma(b_2)}{\Gamma(a_1)}\frac{e^{2z^{\frac{1}{2}}e^{i\frac{\pi}{2}}}(ze^{i\pi})^{\frac{a_1-b_1-b_2+\frac{1}{2}}{2}}}{\sqrt{\pi}}\left\{\sum\limits_{n=0}^{S-1}\frac{\mu_n}{2^{n+1}}(ze^{i\pi})^{-n}+O(|z|^{-S})\right\},\hspace{.12cm}
\label{eq10}
\end{align}
where $a_1$, $b_1$ and $b_2$ are constants and the coefficient $\mu_n$ is given by formula 16.11.4 in \cite{N}.

We then set $z=\frac{\lambda^2 x^2}{4}$, $a_1=\frac{1}{2}$, $b_1=\frac{3}{2}$ and $b_2=\frac{3}{2}$, and we obtain
\begin{eqnarray}
{}_1F_2\left(\frac{1}{2};\frac{3}{2},\frac{3}{2};-\frac{\lambda^2 x^2}{4}\right)=\frac{\pi}{2}\left(\lambda^2 x^2\right)^{-
\frac{1}{2}}\left\{\sum\limits_{n=0}^{R-1}\frac{\left(\frac{1}{2}\right)_n }{n!}\left(i\frac{\lambda x}{2}\right)^{-2n}+O\left(\left|\frac{\lambda x}{2}\right|^{-2R}\right)\right\}\nonumber\\
-\frac{\sqrt{\pi}}{\lambda^2 x^2}\frac{e^{ -i{\lambda x}}}{2}\left\{\sum\limits_{n=0}^{S-1}\frac{\mu_n}{2^n}\left(-i\frac{\lambda x}{2}\right)^{-2n}+O\left(\left|\frac{\lambda x}{2}\right|^{-2S}\right)\right\}
\nonumber\\-\frac{\sqrt{\pi}}{\lambda^2 x^2}\frac{e^{ i{\lambda x}}}{2}\left\{\sum\limits_{n=0}^{S-1}\frac{\mu_n}{2^n}\left(i\frac{\lambda x}{2}\right)^{-2n}+O\left(\left|\frac{\lambda x}{2}\right|^{-2S}\right)\right\}.%\hspace{2cm}
\label{eq11}
\end{eqnarray}
Then, for $|x|\gg1$,
\begin{equation}
\frac{\pi}{2}\left(\lambda^2 x^2\right)^{-
\frac{1}{2}}\left\{\sum\limits_{n=0}^{R-1}\frac{\left(\frac{1}{2}\right)_n }{n!}\left(i\frac{\lambda x}{2}\right)^{-2n}+O\left(\left|\frac{\lambda x}{2}\right|^{-2R}\right)\right\}\sim \frac{\pi}{2\lambda |x|},
\label{eq12}
\end{equation}
while
\begin{multline}
-\frac{\sqrt{\pi}}{\lambda^2 x^2}\frac{e^{ i{\lambda x}}}{2}\left\{\sum\limits_{n=0}^{S-1}\frac{\mu_n}{2^n}\left(-i\frac{\lambda x}{2}\right)^{-2n}
-\sum\limits_{n=0}^{S-1}\frac{\mu_n}{2^n}\left(i\frac{\lambda x}{2}\right)^{-2n}+O\left(\left|\frac{\lambda x}{2}\right|^{-2S}\right)\right\}\\
\sim \frac{\sqrt{\pi}}{\left({\lambda x}\right)^2}\frac{e^{i{\lambda x}}+e^{ -i{\lambda x}}}{2}=
{\sqrt{\pi}}\frac{\cos{\left({\lambda x}\right)}}{\left({\lambda x}\right)^2}.
\label{eq13}
\end{multline}
We then obtain,
\begin{equation}
x{}_1F_2\left(\frac{1}{2};\frac{3}{2},\frac{3}{2};-\frac{\lambda^2 x^2}{4}\right)\sim\frac{\pi}{2\lambda}\frac{x}{|x|}-\frac{\sqrt{\pi}}{\lambda}\frac{\cos{\left({\lambda x}\right)}}{\lambda x}, \hspace{.12cm} |x|\rightarrow\infty.
\label{eq14}
\end{equation}
Hence,
\begin{equation}
G(-\infty)=\lim_{x\rightarrow-\infty}x{}_1F_2\left(\frac{1}{2};\frac{3}{2},\frac{3}{2};-\frac{\lambda^2 x^2}{4}\right)
=\lim_{x\rightarrow-\infty}\left(\frac{\pi}{2\lambda}\frac{x}{|x|}-\frac{\sqrt{\pi}}{\lambda}\frac{\cos{\left({\lambda x}\right)}}{\lambda x}\right)=-\frac{\pi}{2\lambda}
\label{eq15}
\end{equation}
and
\begin{equation}
G(+\infty)=\lim_{x\rightarrow+\infty}x{}_1F_2\left(\frac{1}{2};\frac{3}{2},\frac{3}{2};-\frac{\lambda^2 x^2}{4}\right)
=\lim_{x\rightarrow+\infty}\left(\frac{\pi}{2\lambda}\frac{x}{|x|}-\frac{\sqrt{\pi}}{\lambda}\frac{\cos{\left({\lambda x}\right)}}{\lambda x}\right)=\frac{\pi}{2\lambda}.
\label{eq16}
\end{equation}
\item By the FTC,
\begin{align}
\int\limits_{-\infty}^{+\infty} \frac{\sin{(\lambda x)
}}{\lambda x} dx &=\lim_{y\rightarrow-\infty}\int\limits_{y}^{0} \frac{\sin{(\lambda x)
}}{\lambda x} dx
+\lim_{y\rightarrow+\infty}\int\limits_{0}^{y} \frac{\sin{(\lambda x)
}}{\lambda x}dx=G(+\infty)-G(-\infty)\nonumber\\&=\lim_{y\rightarrow+\infty}y\hspace{.1cm}{}_1F_2\left(\frac{1}{2};\frac{3}{2},\frac{3}{2};-\frac{\lambda^2 y^2}{4}\right)
-\lim_{y\rightarrow-\infty}y\hspace{.1cm}{}_1F_2\left(\frac{1}{2};\frac{3}{2},\frac{3}{2};-\frac{\lambda^2 y^2}{4}\right)
=\frac{\pi}{\lambda}.
\label{eq17}
\end{align}
\end{enumerate}
We now verify whether this is correct or not using Fubini's Theorem \cite{B}. We first observe that
\begin{equation}
\int\limits_{-\infty}^{+\infty} \frac{\sin{(\lambda x)
}}{\lambda x} dx=2\int\limits_{0}^{+\infty} \frac{\sin{(\lambda x),
}}{\lambda x} dx
\label{eq18}
\end{equation}
since the integrand is an even function. We have in terms of double integrals and using Fubini's Theorem that
\begin{equation}
\int\limits_{0}^{+\infty} \frac{\sin{(\lambda x)
}}{\lambda x} dx=\frac{1}{\lambda}\int\limits_{0}^{+\infty}\int\limits_{0}^{+\infty}e^{-s x} {\sin{(\lambda x)}} ds dx=\int\limits_{0}^{+\infty}\int\limits_{0}^{+\infty}e^{-s x} {\sin{(\lambda x)}}  dx ds.
\label{eq19}
\end{equation}
Now using the fact that the inside integral in (\ref{eq19}) is the Laplace transform of $\sin{(\lambda x)}$ \cite{AS} yields
\begin{align}
\int\limits_{0}^{+\infty}\int\limits_{0}^{+\infty}e^{-s x} {\sin{(\lambda x)}}  dx ds=\int\limits_{0}^{+\infty}\frac{\lambda}{s^2+\lambda^2}d s=\arctan{(+\infty)}-\arctan{0}=\frac{\pi}{2}.
\label{eq21}
\end{align}
Hence, \begin{equation}
\int\limits_{-\infty}^{+\infty} \frac{\sin{(\lambda x)
}}{\lambda x} dx=2\int\limits_{0}^{+\infty} \frac{\sin{(\lambda x)
}}{\lambda x} dx=\frac{2}{\lambda}\int\limits_{0}^{+\infty}\int\limits_{0}^{+\infty}e^{-s x} {\sin{(\lambda x)}} ds dx=2\frac{\pi}{2\lambda}=\frac{\pi}{
\lambda}
\label{eq22}
\end{equation}
as before. 
\hfill$\square$\\[1ex]

Setting $\lambda=1$ as in example 1, Lemma \ref{lemen:2} gives $\lim\limits_{x\rightarrow-\infty}G(x)=-\theta=-\pi/2$ while $\lim\limits_{x\rightarrow+\infty}G(x)=\theta=\pi/2$. And these are the exact values of $G(x)$ as $x\rightarrow\pm \infty$ in Figure \ref{fig1}.

\begin{teoen} If $\beta\ge1$ and $\alpha<\beta+1$, then the function
$$G(x)=\frac{x^{\beta-\alpha+1}}{\beta-\alpha+1}\hspace{.075cm} _{1}F_{2}\left(-\frac{\alpha}{2\beta}+\frac{1}{2\beta}+\frac{1}{2};-\frac{\alpha}{2\beta}+\frac{1}{2\beta}+\frac{3}{2},\frac{3}{2};-\frac{\lambda^2 x^{2\beta}}{4}\right),$$
where ${}_1F_2$ is a hypergeometric function \cite{AS} and $\lambda$ is an arbitrarily constant, is the antiderivative of the function $g(x)=\frac{\sin{(\lambda x^\beta)
}}{\lambda x^\alpha}$. Thus,
\begin{align}
\text{Si}_{\beta,\alpha}&=\int\frac{\sin{(\lambda x^\beta)}}{\lambda x^\alpha}dx =\frac{x^{\beta-\alpha+1}}{\beta-\alpha+1}\hspace{.075cm} _{1}F_{2}\left(-\frac{\alpha}{2\beta}+\frac{1}{2\beta}+\frac{1}{2};-\frac{\alpha}{2\beta}+\frac{1}{2\beta}+\frac{3}{2},\frac{3}{2};-\frac{\lambda^2 x^{2\beta}}{4}\right)+C.
\label{eq23}
\end{align}
And for $|x|\gg1$,
\begin{multline}
\frac{x^{\beta-\alpha+1}}{\beta-\alpha+1}\hspace{.075cm}_{1}F_{2}\left(-\frac{\alpha}{2\beta}+\frac{1}{2\beta}+\frac{1}{2};-\frac{\alpha}{2\beta}+\frac{1}{2\beta}+\frac{3}{2},\frac{3}{2};-\frac{\lambda^2 x^{2\beta}}{4}\right)\\\sim\frac{\left(\frac{2}{\lambda}\right)^{1+\frac{1}{\beta}-\frac{\alpha}{\beta}}}{\beta-\alpha+1}\frac{\Gamma\left(-\frac{\alpha}{2\beta}+\frac{1}{2\beta}+\frac{3}{2}\right)}
{\Gamma\left(1+\frac{\alpha}{2\beta}-\frac{1}{2\beta}\right)}\frac{\sqrt{\pi}}{2}
\frac{x^{\beta-\alpha+1}}{|x|^{\beta-\alpha+1}}
-\frac{\beta-\alpha+1}{\beta}\frac{\sqrt{\pi}}{\lambda^2}\frac{\cos{\left({\lambda x^\beta}\right)}}{x^{\beta+\alpha-1}}. 
\label{eq25}
\end{multline}
\label{teoen:1}
\end{teoen}

\proofen
\begin{align}
\text{Si}_{\beta,\alpha}&=\int g(x) dx=\int \frac{\sin{(\lambda x^\beta)}}{\lambda x^\alpha}dx
=\int\frac{1}{\lambda x^\alpha}\sum\limits_{n=0}^{\infty}(-1)^n\frac{(\lambda x^\beta)^{2n+1}}{(2n+1)!}dx
\nonumber\\ &=\sum\limits_{n=0}^{\infty}(-1)^n\frac{\lambda^{2n}}{(2n+1)!}\int{x^{2\beta n+\beta-\alpha}}dx
=\lambda\sum\limits_{n=0}^{\infty}(-1)^n\frac{\lambda^{2n}}{(2n+1)!}\int x^{2\beta n+\beta -\alpha} dx
\nonumber\\&=\sum\limits_{n=0}^{\infty}(-1)^n\frac{\lambda^{2n}}{(2n+1)!}\frac{x^{2\beta n+\beta-\alpha+1}}{2\beta n+\beta-\alpha+1}+C
\nonumber\\&=\frac{x^{\beta-\alpha+1}}{2\beta}\sum\limits_{n=0}^{\infty}(-1)^n\frac{\lambda^{2n}}{(2n+1)!}\frac{x^{2\beta n}}{n-\frac{\alpha}{2\beta}+\frac{1}{2\beta}+\frac{1}{2}}+C
\nonumber\\&=\frac{x^{\beta-\alpha+1}}{2\beta}\sum\limits_{n=0}^{\infty}\frac{\Gamma\left(n-\frac{\alpha}{2\beta}+\frac{1}{2\beta}+\frac{1}{2}\right)}{\Gamma(2n+2)\Gamma\left(n-\frac{\alpha}{2\beta}+\frac{1}{2\beta}+\frac{3}{2}\right)}(-\lambda^2 x^{2\beta})^{n}+C
\nonumber\\&=\frac{x^{\beta-\alpha+1}}{\beta-\alpha+1}\sum\limits_{n=0}^{\infty}\frac{\left(-\frac{\alpha}{2\beta}+\frac{1}{2\beta}+\frac{1}{2}\right)_n}
{\left(\frac{3}{2}\right)_n\left(-\frac{\alpha}{2\beta}+\frac{1}{2\beta}+\frac{3}{2}\right)_n}\frac{\left(-\frac{\lambda^2 x^{2\beta}}{4}\right)^{n}}{n!}+C
\nonumber\\&=\frac{x^{\beta-\alpha+1}}{\beta-\alpha+1}\hspace{.075cm} _{1}F_{2}\left(-\frac{\alpha}{2\beta}+\frac{1}{2\beta}+\frac{1}{2};-\frac{\alpha}{2\beta}+\frac{1}{2\beta}+\frac{3}{2},\frac{3}{2};-\frac{\lambda^2 x^{2\beta}}{4}\right)+C
=G(x)+C.
\label{eq27}
\end{align}
To prove (\ref{eq25}), we use the asymptotic formula for the hypergeometric function $_{1}F_{2}$ in equation (\ref{eq10}), and proceed as in Lemma \ref{lemen:2}. 
\hfill$\square$\\[1ex]

Beside, we can show as above that if $\beta\ge1$ and $\alpha<\beta+1$, then
\begin{align}
\int \frac{\sinh{(\lambda x^\beta)}}{\lambda x^\alpha}dx
=\frac{x^{\beta-\alpha+1}}{\beta-\alpha+1}\hspace{.075cm} _{1}F_{2}\left(-\frac{\alpha}{2\beta}+\frac{1}{2\beta}+\frac{1}{2};-\frac{\alpha}{2\beta}+\frac{1}{2\beta}+\frac{3}{2},\frac{3}{2};\frac{\lambda^2 x^{2\beta}}{4}\right)+C.
\label{eq28}
\end{align}

\begin{coren}
Let $\beta=\alpha$. If $\alpha\ge1$, then
\begin{equation}
\int\limits_{-\infty}^{0} \frac{\sin{(\lambda x^\alpha)
}}{\lambda x^\alpha}dx=G(0)-G(-\infty)=\left(\frac{2}{\lambda}\right)^{\frac{1}{\alpha}}\frac{\Gamma\left(\frac{1}{2\alpha}+1\right)}{\Gamma\left(\frac{3}{2}-\frac{1}{2\alpha}\right)}\frac{\sqrt{\pi}}{2},
\label{eq29}
\end{equation}
\begin{equation}
\int\limits_{0}^{+\infty} \frac{\sin{(\lambda x^\alpha)
}}{\lambda x^\alpha}dx=G(+\infty)-G(0)=\left(\frac{2}{\lambda}\right)^{\frac{1}{\alpha}}\frac{\Gamma\left(\frac{1}{2\alpha}+1\right)}{\Gamma\left(\frac{3}{2}-\frac{1}{2\alpha}\right)}\frac{\sqrt{\pi}}{2}
\label{eq30}
\end{equation}
and
\begin{equation}
\int\limits_{-\infty}^{+\infty} \frac{\sin{(\lambda x^\alpha)
}}{\lambda x^\alpha}dx=G(+\infty)-G(-\infty)=\left(\frac{2}{\lambda}\right)^{\frac{1}{\alpha}}\frac{\Gamma\left(\frac{1}{2\alpha}+1\right)}{\Gamma\left(\frac{3}{2}-\frac{1}{2\alpha}\right)}\sqrt{\pi}.
\label{eq31}
\end{equation}
\label{coren:1}
\end{coren}

\proofen
If $\beta=\alpha$, Theorem \ref{teoen:1} gives
\begin{align}
G(-\infty)&=\lim_{x\rightarrow-\infty}x{}_1F_2\left(\frac{1}{2\alpha};\frac{1}{2\alpha}+1,\frac{3}{2};-\frac{\lambda^2 x^{2\alpha}}{4}\right)
\nonumber\\&=\lim_{x\rightarrow-\infty}\left(\left(\frac{2}{\lambda}\right)^{\frac{1}{\alpha}}\frac{\sqrt{\pi}}{2}\frac{\Gamma\left(\frac{1}{2\alpha}+1\right)}{\Gamma\left(\frac{3}{2}-\frac{1}{2\alpha}\right)}\frac{x}{|x|}
-\frac{\sqrt{\pi}}{\alpha\lambda^2}\frac{\cos{\left({\lambda x^\alpha}\right)}}{x^{2\alpha-1}}\right)
=-\left(\frac{2}{\lambda}\right)^{\frac{1}{\alpha}}\frac{\Gamma\left(\frac{1}{2\alpha}+1\right)}{\Gamma\left(\frac{3}{2}-\frac{1}{2\alpha}\right)}\frac{\sqrt{\pi}}{2}
\label{eq32}
\end{align}
and
\begin{align}
G(+\infty)&=\lim_{x\rightarrow+\infty}x{}_1F_2\left(\frac{1}{2\alpha};\frac{1}{2\alpha}+1,\frac{3}{2};-\frac{\lambda^2 x^{2\alpha}}{4}\right)
\nonumber\\&=\lim_{x\rightarrow+\infty}\left(\left(\frac{2}{\lambda}\right)^{\frac{1}{\alpha}}\frac{\sqrt{\pi}}{2}\frac{\Gamma\left(\frac{1}{2\alpha}+1\right)}{\Gamma\left(\frac{3}{2}-\frac{1}{2\alpha}\right)}\frac{x}{|x|}
-\frac{\sqrt{\pi}}{\alpha\lambda^2}\frac{\cos{\left({\lambda x^\alpha}\right)}}{x^{2\alpha-1}}\right)
=\left(\frac{2}{\lambda}\right)^{\frac{1}{\alpha}}\frac{\Gamma\left(\frac{1}{2\alpha}+1\right)}{\Gamma\left(\frac{3}{2}-\frac{1}{2\alpha}\right)}\frac{\sqrt{\pi}}{2}.
\label{eq33}
\end{align}
Hence, by the FTC,
\begin{equation}
\int\limits_{-\infty}^{0} \frac{\sin{(\lambda x^\alpha)
}}{\lambda x^\alpha}dx=G(0)-G(-\infty)=0-\left(-\left(\frac{2}{\lambda}\right)^{\frac{1}{\alpha}}\frac{\Gamma\left(\frac{1}{2\alpha}+1\right)}{\Gamma\left(\frac{3}{2}-\frac{1}{2\alpha}\right)}\frac{\sqrt{\pi}}{2}\right)
=\left(\frac{2}{\lambda}\right)^{\frac{1}{\alpha}}\frac{\Gamma\left(\frac{1}{2\alpha}+1\right)}{\Gamma\left(\frac{3}{2}-\frac{1}{2\alpha}\right)}\frac{\sqrt{\pi}}{2},
\label{eq34}
\end{equation}
\begin{equation}
\int\limits_{0}^{+\infty} \frac{\sin{(\lambda x^\alpha)
}}{\lambda x^\alpha}dx=G(+\infty)-G(0)=\left(\frac{2}{\lambda}\right)^{\frac{1}{\alpha}}\frac{\Gamma\left(\frac{1}{2\alpha}+1\right)}{\Gamma\left(\frac{3}{2}-\frac{1}{2\alpha}\right)}
\frac{\sqrt{\pi}}{2}-0=\left(\frac{2}{\lambda}\right)^{\frac{1}{\alpha}}\frac{\Gamma\left(\frac{1}{2\alpha}+1\right)}{\Gamma\left(\frac{3}{2}-\frac{1}{2\alpha}\right)}\frac{\sqrt{\pi}}{2}.
\label{eq35}
\end{equation}
And combining (\ref{eq34}) and (\ref{eq35}) gives (\ref{eq31}).
\hfill $\square$

\begin{teoen}
If $\beta \ge1$ and $\alpha <\beta+1$, then the FTC gives
\begin{equation}
\int\limits_{A}^{B} \frac{\sin{(\lambda x^\beta)
}}{\lambda x^\alpha}dx=G(B)-G(A),
\label{eq36}
\end{equation}
for any $A$ and any $B$, and where $G$ is given in Theorem \ref{teoen:1}.
\label{teoen:2}
\end{teoen}

\proofen
Equation (\ref{eq36}) holds by Theorem \ref{teoen:1}, Corollary \ref{coren:1} and Lemma \ref{lemen:2}. Since the FTC works for $A=-\infty$ and $B=0$ in (\ref{eq29}), $A=0$ and $B=+\infty$ in (\ref{eq30}) and  $A=-\infty$ and $B=+\infty$ in (\ref{eq31}) by Corollary \ref{coren:1} for any $\beta=\alpha \ge1$ and by Lemma \ref{lemen:2} for $\beta=\alpha=1$, then it has to work for other values of $A,B\in \mathbb{R}$ and for $\beta\ge1$ and $\alpha<\beta+1$ since the case with $\beta=\alpha\ge1$ is derived from the case with $\beta\ge1$ and $\alpha<\beta+1$.
\hfill $\square$

\begin{teoen}
Let $\beta\ge1$, then the function
$G(x)=\ln|x|-\frac{\left(\frac{\lambda x^{\beta}}{2}\right)^2}{6\beta}\ _{2}F_{3}\left(1,1;2,2,\frac{5}{2};-\frac{\lambda^2 x^{2\beta}}{4}\right),$
where ${}_2F_3$ is a hypergeometric function \cite{AS} and $\lambda$ is an arbitrarily constant, is the antiderivative of the function $g(x)=\frac{\sin{(\lambda x^\beta)
}}{\lambda x^{\beta+1}}$. Thus,
\begin{align}
\text{Si}_{\beta,\beta+1}&=\int\frac{\sin{(\lambda x^\beta)}}{\lambda x^{\beta+1}}dx =\ln|x|-\frac{\left(\frac{\lambda x^{\beta}}{2}\right)^2}{6\beta}\ _{2}F_{3}\left(1,1;2,2,\frac{5}{2};-\frac{\lambda^2 x^{2\beta}}{4}\right)+C.
\label{eq37}
\end{align}
\label{teoen:3}
\end{teoen}

\proofen
\begin{align}
\text{Si}_{\beta,\beta+1}&=\int g(x) dx=\int \frac{\sin{(\lambda x^\beta)}}{\lambda x^{\beta+1}}dx
=\int\frac{1}{\lambda x^{\beta+1}}\sum\limits_{n=0}^{\infty}(-1)^n\frac{(\lambda x^\beta)^{2n+1}}{(2n+1)!}dx
\nonumber\\ &=\sum\limits_{n=0}^{\infty}(-1)^n\frac{\lambda^{2n}}{(2n+1)!}\int{x^{2\beta n-1}}dx
=\int \frac{dx}{x}+\sum\limits_{n=1}^{\infty}(-1)^n\frac{\lambda^{2n}}{(2n+1)!}\int x^{2\beta n-1} dx
\nonumber\\&=\ln|x|+\sum\limits_{n=0}^{\infty}(-1)^{n+1}\frac{\lambda^{2n+2}}{(2n+3)!}\frac{x^{2\beta n+2\beta}}{2\beta n+2\beta}+C
\nonumber\\&=\ln|x|-\frac{\lambda^2 x^{2\beta}}{2\beta}\sum\limits_{n=0}^{\infty}(-1)^n\frac{\lambda^{2n}}{(2n+3)!}\frac{x^{2\beta n}}{n+1}+C
\nonumber\\&=\ln|x|-\frac{\lambda^2x^{2\beta}}{2\beta}\sum\limits_{n=0}^{\infty}\frac{\left(\Gamma\left(n+1\right)\right)^2}{\Gamma(2n+4)\Gamma\left(n+2\right)}\frac{\left(-\lambda^2 x^{2\beta}\right)^{n}}{n!}+C
\nonumber\\&=\ln|x|-\frac{\left(\frac{\lambda x^{\beta}}{2}\right)^2}{6\beta}\sum\limits_{n=0}^{\infty}\frac{(1)_n(1)_n}
{(2)_n (2)_n\left(\frac{5}{2}\right)_n}\frac{\left(-\frac{\lambda^2 x^{2\beta}}{4}\right)^{n}}{n!}+C
\nonumber\\&=\ln|x|-\frac{\left(\frac{\lambda x^{\beta}}{2}\right)^2}{6\beta}\ _{2}F_{3}\left(1,1;2,2,\frac{5}{2};-\frac{\lambda^2 x^{2\beta}}{4}\right)+C
=G(x)+C.
\label{eq38}
\end{align} \hfill $\square$

\section{Evaluation of the cosine integral and related integrals}\label{sec:3}
\begin{teoen}
If $\beta\ge1$ and $\alpha <2\beta+1$, then the function
$$G(x)=\frac{1}{\lambda}\frac{x^{1-\alpha}}{1-\alpha}-\frac{\lambda x^{2\beta-\alpha+1}}{2\beta-\alpha+1}\hspace{.075cm} _{2}F_{3}\left(1,-\frac{\alpha}{2\beta}+\frac{1}{2\beta}+1;-\frac{\alpha}{2\beta}+\frac{1}{2\beta}+2,\frac{3}{2},2;-\frac{\lambda^2 x^{2\beta}}{4}\right),$$
where ${}_2F_3$ is a hypergeometric function \cite{AS} and $\lambda$ is an arbitrarily constant, is the antiderivative of the function $g(x)=\frac{\cos{(\lambda x^\beta)
}}{\lambda x^\alpha}$. Thus,
\begin{multline}
\int\frac{\cos{(\lambda x^\beta)}}{\lambda x^\alpha}dx=\frac{1}{\lambda}\frac{x^{1-\alpha}}{1-\alpha}-\frac{1}{2}\frac{\lambda x^{2\beta-\alpha+1}}{2\beta-\alpha+1}\hspace{.075cm} _{2}F_{3}\left(1,-\frac{\alpha}{2\beta}+\frac{1}{2\beta}+1;-\frac{\alpha}{2\beta}+\frac{1}{2\beta}+2,\frac{3}{2},2;-\frac{\lambda^2 x^{2\beta}}{4}\right)+C,
\label{eq39}
\end{multline}
and for $|x|\gg1$,
\begin{multline}
\frac{\lambda x^{2\beta-\alpha+1}}{2\beta-\alpha+1}\hspace{.075cm} _{2}F_{3}\left(1,-\frac{\alpha}{2\beta}+\frac{1}{2\beta}+1;-\frac{\alpha}{2\beta}+\frac{1}{2\beta}+2,\frac{3}{2},2;-\frac{\lambda^2 x^{2\beta}}{4}\right)\\\sim
\frac{\sqrt{\pi}\lambda}{2\beta}\Gamma\left(-\frac{\alpha}{2\beta}+\frac{1}{2\beta}+1\right)\left(\frac{2}{\lambda }\right)^{-\frac{\alpha}{\beta}+\frac{1}{\beta}+2}+\frac{\sqrt{\pi}}{\lambda\beta}x^{-\alpha+1}
+\frac{2}{\lambda^2\beta}\frac{\cos(\lambda x^{\beta})}{x^{\beta+\alpha-1}}.
\label{eq40}
\end{multline}
\label{teoen:4}
\end{teoen}

\proofen
If $\beta\ge1$ and $\alpha<2\beta+1$,
\begin{align}
\int g(x) dx&=\int \frac{\cos{(\lambda x^\beta)
}}{\lambda x^\alpha}dx =\int\frac{1}{\lambda x^\alpha}\sum\limits_{n=0}^{\infty}(-1)^n\frac{(\lambda x^\beta)^{2n}}{(2n)!}dx
\nonumber\\ &=\int\frac{1}{\lambda x^\alpha}dx+\frac{1}{\lambda}\int\sum\limits_{n=1}^{\infty}(-1)^n\frac{\lambda^{2n}}{(2n)!}x^{2\beta n-\alpha}dx
\nonumber\\&=\frac{1}{\lambda}\frac{x^{1-\alpha}}{1-\alpha}-\frac{1}{\lambda}\sum\limits_{n=0}^{\infty}(-1)^n\frac{\lambda^{2n+2}}{(2n+2)!}\int x^{2\beta n+2\beta-\alpha} dx
\nonumber\\&=\frac{1}{\lambda}\frac{x^{1-\alpha}}{1-\alpha}-\lambda\sum\limits_{n=0}^{\infty}(-1)^n\frac{\lambda^{2n}}{(2n+2)!}\frac{x^{2\beta n+2\beta-\alpha+1}}{2\beta n+2\beta-\alpha+1}+C
\nonumber\\&=\frac{1}{\lambda}\frac{x^{1-\alpha}}{1-\alpha}-\frac{\lambda x^{2\beta-\alpha+1}}{2\beta}\sum\limits_{n=0}^{\infty}
\frac{\Gamma\left(n-\frac{\alpha}{2\beta}+\frac{1}{2\beta}+1\right)}{\Gamma(2n+3)\Gamma\left(n-\frac{\alpha}{2\beta}+\frac{1}{2\beta}+2\right)}(-\lambda^2 x^{2\beta})^{n}+C
\nonumber\\&=\frac{1}{\lambda}\frac{x^{1-\alpha}}{1-\alpha}-\frac{1}{2}\frac{\lambda x^{2\beta-\alpha+1}}{2\beta-\alpha+1}-\sum\limits_{n=0}^{\infty}\frac{(1)_n\left(-\frac{\alpha}{2\beta}+\frac{1}{2\beta}+1\right)_n}{\left(\frac{3}{2}\right)_n (2)_n\left(-\frac{\alpha}{2\beta}+\frac{1}{2\beta}+2\right)_n}\frac{\left(-\frac{\lambda^2 x^{2\beta}}{4}\right)^{n}}{n!}+C
\nonumber\\&=\frac{1}{\lambda}\frac{x^{1-\alpha}}{1-\alpha}-\frac{1}{2}\frac{\lambda x^{2\beta-\alpha+1}}{2\beta-\alpha+1}\hspace{.075cm} _{2}F_{3}\left(1,-\frac{\alpha}{2\beta}+\frac{1}{2\beta}+1;-\frac{\alpha}{2\beta}+\frac{1}{2\beta}+2,\frac{3}{2},2;-\frac{\lambda^2 x^{2\beta}}{4}\right)+C.
%\nonumber\\&=G(x)+C.
\label{eq41}
\end{align}
To prove (\ref{eq40}), we use the asymptotic expression of ${}_2F_3\left(a_1,a_2;b_1,b_2,b_3;-z\right)$ for $|z|\gg 1$, where $a_1,a_2,b_1,b_2$ and $b_3$ are constants, and $-\pi\le\mbox{arg}\hspace{.1cm} z\le\pi$ . It can be obtained using formulas 16.11.1, 16.11.2 and 16.11.8 in \cite{N} and is given by
\begin{align}
&{}_2F_3\left(a_1,a_2;b_1,b_2,b_3;-z\right)=
\nonumber \\& \frac{\Gamma(b_1)\Gamma(b_2)\Gamma(b_3)}{\Gamma(a_2)}z^{-a_1} \left\{\sum\limits_{n=0}^{R-1}\frac{(a_1)_n\Gamma(a_1-a_2-n) }{\Gamma(b_1-a_1-n)\Gamma(b_2-a_1-n)\Gamma(b_3-a_1-n)}\frac{(-z)^{-n}}{n!}+O(|z|^{-R})\right\}\nonumber\\&
+\frac{\Gamma(b_1)\Gamma(b_2)\Gamma(b_3)}{\Gamma(a_1)}z^{-a_2}\left\{\sum\limits_{n=0}^{R-1}\frac{(a_2)_n\Gamma(a_2-a_1-n) }{\Gamma(b_1-a_2-n)\Gamma(b_2-a_2-n)\Gamma(b_3-a_2-n)}\frac{(-z)^{-n}}{n!}+O(|z|^{-R})\right\}\nonumber\\&
+\frac{\Gamma(b_1)\Gamma(b_2)\Gamma(b_3)}{\Gamma(a_1)\Gamma(a_2)}\frac{e^{2z^{\frac{1}{2}}e^{-i\frac{\pi}{2}}}(ze^{-i\pi})^{\frac{a_1+a_2-b_1-b_2-b_3+\frac{1}{2}}{2}}}{\sqrt{\pi}}\left\{\sum\limits_{n=0}^{S-1}\frac{\mu_n}{2^{n+1}}(ze^{-i\pi})^{-n}+O(|z|^{-S})\right\}
\nonumber\\&+\frac{\Gamma(b_1)\Gamma(b_2)\Gamma(b_3)}{\Gamma(a_1)\Gamma(a_2)} \frac{e^{2z^{\frac{1}{2}}e^{i\frac{\pi}{2}}}(ze^{i\pi})^{\frac{a_1+a_2-b_1-b_2-b_3+\frac{1}{2}}{2}}}{\sqrt{\pi}}\left\{\sum\limits_{n=0}^{S-1}\frac{\mu_n}{2^{n+1}}(ze^{i\pi})^{-n}+O(|z|^{-S})\right\},
\label{eq42}
\end{align}
where the coefficient $\mu_n$ is given by formula 16.11.4 in \cite{N}.

We now set $z=\frac{\lambda^2 x^{2\beta}}{4}$, $a_1=1$, $a_2=-\frac{\alpha}{2\beta}+\frac{1}{2\beta}+1$,
$b_1=-\frac{\alpha}{2\beta}+\frac{1}{2\beta}+2$, $b_2=\frac{3}{2}$ and $b_3={2}$ in (\ref{eq42}) to obtain
\begin{align}
& _{2}F_{3}\left(1,-\frac{\alpha}{2\beta}+\frac{1}{2\beta}+1;-\frac{\alpha}{2\beta}+\frac{1}{2\beta}+2,\frac{3}{2},2;-\frac{\lambda^2 x^{2\beta}}{4}\right)\nonumber \\ &\sim \frac{\sqrt{\pi}}{\lambda^2}\left(-\frac{\alpha}{\beta}+\frac{1}{\beta}+{2}\right)\frac{1}{x^{2\beta}}
+\frac{\sqrt{\pi}}{2}\Gamma\left(-\frac{\alpha}{2\beta}+\frac{1}{2\beta}+2\right)\left(\frac{2}{\lambda x^{\beta}}\right)^{-\frac{\alpha}{\beta}+\frac{1}{\beta}+2}
+\frac{2}{\lambda^3}\left(-\frac{\alpha}{\beta}+\frac{1}{\beta}+{2}\right)\frac{\cos(\lambda x^{\beta})}{x^{3\beta}}.
\label{eq43}
\end{align}
Hence, multiplying (\ref{eq43}) with $\frac{\lambda x^{2\beta-\alpha+1}}{2\beta-\alpha+1}$ gives (\ref{eq40}). 
\hfill $\square$

On the other hand, we can show as above that if $\beta\ge1$ and $\alpha<2\beta+1$, then
\begin{equation}
\int\frac{\cosh{(\lambda x^\beta)}}{\lambda x^\alpha}dx=\frac{1}{\lambda}\frac{x^{1-\alpha}}{1-\alpha}+\frac{1}{2}\frac{\lambda x^{2\beta-\alpha+1}}{2\beta-\alpha+1}\hspace{.075cm} _{2}F_{3}\left(1,-\frac{\alpha}{2\beta}+\frac{1}{2\beta}+1;-\frac{\alpha}{2\beta}+\frac{1}{2\beta}+2,\frac{3}{2},2;\frac{\lambda^2 x^{2\beta}}{4}\right)+C.
\label{eq44}
\end{equation}

\begin{teoen}
Let $\beta\ge1$, then the function $G(x)=-\frac{1}{2\lambda\beta x^{2\beta}}-\frac{\lambda}{2}\ln|x|+\frac{\lambda}{6\beta}\left(\frac{\lambda x^{\beta}}{4}\right)^2\hspace{.075cm} _{2}F_{3}\left(1,1;2,\frac{5}{2},3;-\frac{\lambda^2 x^{2\beta}}{4}\right)$, where ${}_2F_3$ is a hypergeometric function \cite{AS} and $\lambda$ is an arbitrarily constant, is the antiderivative of the function $g(x)=\frac{\cos{(\lambda x^\beta)
}}{\lambda x^{2\beta+1}}$. Thus,
\begin{equation}
\text{Ci}_{\beta,2\beta+1}=\int\frac{\cos{(\lambda x^\beta)}}{\lambda x^{2\beta+1}}dx=-\frac{1}{2\lambda\beta x^{2\beta}}-\frac{\lambda}{2}\ln|x|+\frac{\lambda}{6\beta}\left(\frac{\lambda x^{\beta}}{4}\right)^2\hspace{.075cm} _{2}F_{3}\left(1,1;2,\frac{5}{2},3;-\frac{\lambda^2 x^{2\beta}}{4}\right)+C.
\label{eq45}
\end{equation}
We also have,
\begin{equation}
\text{Ci}_{\beta,1}=\int\frac{\cos{(\lambda x^\beta)}}{\lambda x}dx=\frac{1}{\lambda}\ln|x|-\frac{\lambda x^{2\beta}}{4\beta}\hspace{.075cm} _{2}F_{3}\left(1,1;\frac{3}{2},2,2;-\frac{\lambda^2 x^{2\beta}}{4}\right)+C.
\label{eq46}
\end{equation}
\label{teoen:5}
\end{teoen}

\proofen
\begin{align}
\text{Ci}_{\beta,2\beta+1}=\int g(x) dx&=\int \frac{\cos{(\lambda x^\beta)
}}{\lambda x^{2\beta+1}}dx=\int\frac{1}{\lambda x^{2\beta+1}}\sum\limits_{n=0}^{\infty}(-1)^n\frac{(\lambda x^\beta)^{2n}}{(2n)!}dx
\nonumber\\ &=\int\frac{1}{\lambda x^{2\beta+1}}dx+\frac{1}{\lambda}\int\sum\limits_{n=1}^{\infty}(-1)^n\frac{\lambda^{2n}}{(2n)!}x^{2\beta n-2\beta-1}dx
\nonumber\\&=-\frac{1}{2\lambda\beta x^{2\beta}}-\frac{1}{\lambda}\sum\limits_{n=0}^{\infty}(-1)^n\frac{\lambda^{2n+2}}{(2n+2)!}\int x^{2\beta n-1} dx
\nonumber\\&=-\frac{1}{2\lambda\beta x^{2\beta}}-\frac{\lambda}{2}\int\frac{dx}{x}-{\lambda}\sum\limits_{n=1}^{\infty}(-1)^n\frac{\lambda^{2n}}{(2n+2)!}\int x^{2\beta n-1} dx
\nonumber\\&=-\frac{1}{2\lambda\beta x^{2\beta}}-\frac{\lambda}{2}\ln|x|+{\lambda^3}\sum\limits_{n=0}^{\infty}(-1)^n\frac{\lambda^{2n}}{(2n+4)!}\int x^{2\beta n+2\beta-1} dx
\nonumber\\&=-\frac{1}{2\lambda\beta x^{2\beta}}-\frac{\lambda}{2}\ln|x|+\lambda^3\sum\limits_{n=0}^{\infty}(-1)^n\frac{\lambda^{2n}}{(2n+4)!}\frac{x^{2\beta n+2\beta}}{2\beta n+2\beta}+C
\nonumber\\&=-\frac{1}{2\lambda\beta x^{2\beta}}-\frac{\lambda}{2}\ln|x|-\frac{\lambda^3 x^{2\beta}}{2\beta}\sum\limits_{n=0}^{\infty}
\frac{(\Gamma\left(n+1\right))^2}{\Gamma(2n+5)\Gamma\left(n+2\right)}\frac{(-\lambda^2 x^{2\beta})^{n}}{n!}+C
\nonumber\\&=-\frac{1}{2\lambda\beta x^{2\beta}}-\frac{\lambda}{2}\ln|x|+\frac{\lambda}{6\beta}\left(\frac{\lambda x^{\beta}}{4}\right)^2\sum\limits_{n=0}^{\infty}\frac{(1)_n\left(1\right)_n}{(2)_n\left(\frac{5}{2}\right)_n\left(3\right)_n}\frac{\left(-\frac{\lambda^2 x^{2\beta}}{4}\right)^{n}}{n!}+C
\nonumber\\&=-\frac{1}{2\lambda\beta x^{2\beta}}-\frac{\lambda}{2}\ln|x|+\frac{\lambda}{6\beta}\left(\frac{\lambda x^{\beta}}{4}\right)^2\hspace{.075cm} _{2}F_{3}\left(1,1;2,\frac{5}{2},3;-\frac{\lambda^2 x^{2\beta}}{4}\right)+C
\label{eq47}
\end{align}
The proof of (\ref{eq46}) is similar, we do not show it here.
\hfill $\square$

\section{Evaluation of some integrals involving $\text{Si}_{\alpha,\beta}$ and $\text{Ci}_{\alpha,\beta}$}\label{sec:4}
The integral $\int\frac{\cos^n{(\lambda x^\beta)}}{\lambda x^\alpha}dx$, where $n$ is a positive integer and $\beta\ge1,\alpha<2\beta+1$, can be written in terms of (\ref{eq39}) and then evaluated.
%\emph{Example~1.} Text...\\[1ex] % Example

\emph{Example~2.}
In this example, the integral $\int\frac{\cos^4{(\lambda x^\beta)}}{\lambda x^\alpha}dx$ is evaluated by linearizing $\cos^4{(\lambda x^\beta)}$.
\begin{align}
&\int\frac{\cos^4{(\lambda x^\beta)}}{\lambda x^\alpha}dx=\frac{1}{8}\int\frac{\cos{(4\lambda x^\beta)}}{\lambda x^\alpha}dx+\frac{1}{2}\int\frac{\cos{(2\lambda x^\beta)}}{\lambda x^\alpha}dx+\frac{3}{8}\int dx=
\nonumber\\ &\frac{1}{8\lambda}\frac{x^{1-\alpha}}{1-\alpha}-\frac{1}{4}\frac{\lambda x^{2\beta-\alpha+1}}{2\beta-\alpha+1}\hspace{.075cm} _{2}F_{3}\left(1,-\frac{\alpha}{2\beta}+\frac{1}{2\beta}+1;-\frac{\alpha}{2\beta}+\frac{1}{2\beta}+2,\frac{3}{2},2;-4\lambda^2 x^{2\beta}\right)
\nonumber\\ &+\frac{1}{2\lambda}\frac{x^{1-\alpha}}{1-\alpha}-\frac{1}{2}\frac{\lambda x^{2\beta-\alpha+1}}{2\beta-\alpha+1}\hspace{.075cm} _{2}F_{3}\left(1,-\frac{\alpha}{2\beta}+\frac{1}{2\beta}+1;-\frac{\alpha}{2\beta}+\frac{1}{2\beta}+2,\frac{3}{2},2;-\lambda^2 x^{2\beta}\right)\nonumber\\ &+\frac{3x}{8}+C.
\label{eq48}
\end{align}
%\label{ex:2}
%\end{exaen}
\\[2ex]
If $\beta\ge1$ and $\alpha<\beta+1$, the integral $\int\frac{\sin^n{(\lambda x^\beta)}}{\lambda x^\alpha}dx$, where $n$ is a positive integer, can be written either in terms of (\ref{eq23}) if $n$ odd, and then evaluated.\\

%\begin{exaen}
\emph{Example~3.}
In this example, the integral $\int\frac{\sin^3{(\lambda x^\beta)}}{\lambda x^\alpha}dx$ is evaluated by linearizing $\sin^3{(\lambda x^\beta)}$. 
\begin{align}
\int\frac{\sin^3{(\lambda x^\beta)}}{\lambda x^\alpha}dx&=-\frac{1}{4}\int\frac{\sin{(3\lambda x^\beta)}}{\lambda x^\alpha}dx+\frac{3}{4}\int\frac{\sin{(\lambda x^\beta)}}{\lambda x^\alpha}dx
\nonumber\\ &=-\frac{1}{4}\frac{x^{\beta-\alpha+1}}{\beta-\alpha+1}\hspace{.075cm} _{1}F_{2}\left(-\frac{\alpha}{2\beta}+\frac{1}{2\beta}+\frac{1}{2};-\frac{\alpha}{2\beta}+\frac{1}{2\beta}+\frac{3}{2},\frac{3}{2};-\frac{9\lambda^2 x^{2\beta}}{4}\right)
\nonumber\\ &\hspace{.4cm}+\frac{3}{4}\frac{x^{\beta-\alpha+1}}{\beta-\alpha+1}\hspace{.075cm} _{1}F_{2}\left(-\frac{\alpha}{2\beta}+\frac{1}{2\beta}+\frac{1}{2};-\frac{\alpha}{2\beta}+\frac{1}{2\beta}+\frac{3}{2},\frac{3}{2};-\frac{\lambda^2 x^{2\beta}}{4}\right)+C.
\label{eq49}
\end{align}
%\label{ex:3}
%\end{exaen}
\\[3ex]

\emph{Example~4.} Let us now evaluate the integrals $\int\sin{\left(\frac{\lambda} {x^\mu}\right)}dx$ and $\int\cos{\left(\frac{\lambda} {x^\mu}\right)}dx$.
\begin{enumerate}
\item The integral $\int\sin{\left(\frac{\lambda} {x^\mu}\right)}dx$ is evaluated using the substitution $u=1/x$ and Theorem \ref{teoen:1} if $\mu>1$. Then, we have 
\begin{align}
\int\sin{\left(\frac{\lambda} {x^\mu}\right)}dx&=-\int\frac{\sin{\left(\lambda u^\mu\right)}}{u^2}du
=-\frac{\lambda u^{\mu-1}}{\mu-1}\hspace{.075cm} _{1}F_{2}\left(-\frac{1}{2\mu}+\frac{1}{2};-\frac{1}{2\mu}+\frac{3}{2},\frac{3}{2};-\frac{\lambda^2 u^{2\mu}}{4}\right)
\nonumber\\ &=-\frac{\lambda\left(\frac{1}{x}\right)^{\mu-1}}{\mu-1}\hspace{.075cm} _{1}F_{2}\left(-\frac{1}{2\mu}+\frac{1}{2};-\frac{1}{2\mu}+\frac{3}{2},\frac{3}{2};-\frac{\lambda^2 }{4 x^{2\mu}}\right)+C, \mu>1.
\label{eq50}
\end{align}
The integral $\int\sin{\left(\frac{\lambda} {x^\mu}\right)}dx$ is evaluated using the substitution $u=1/x$ and Theorem \ref{teoen:3} if $\mu=1$. Then, we have
\begin{align}
\int\sin{\left(\frac{\lambda} {x}\right)}dx&=-\int\frac{\sin{\left(\lambda u\right)}}{u^2}du
=-\ln|u|+\frac{\left(\frac{\lambda u}{2}\right)^2}{6}\ _{2}F_{3}\left(1,1;2,2,\frac{5}{2};-\frac{\lambda^2 u^{2}}{4}\right)
\nonumber\\ &=\ln|x|+\frac{\left(\frac{\lambda}{2 x}\right)^2}{6}\ _{2}F_{3}\left(1,1;2,2,\frac{5}{2};-\frac{\lambda^2 }{4 x^{2}}\right)+C.
\label{eq50-1}
\end{align}
\item Making the substitution $u=1/x$ and applying Theorem \ref{teoen:4} gives 
\begin{align}
\int\cos{\left(\frac{\lambda} {x^\mu}\right)}dx&=-\int\frac{\cos{\left(\lambda u^\mu\right)}}{u^2}du
=\frac{1}{u}+\frac{\lambda u^{2\mu-1}}{2\mu-1}\hspace{.075cm} _{2}F_{3}\left(1,-\frac{1}{2\mu}+1;-\frac{1}{2\mu}+2,\frac{3}{2},2;-\frac{\lambda^2 u^{2\mu}}{4}\right)
\nonumber\\ &=x+\frac{\lambda\left(\frac{1}{x}\right)^{2\mu-1}}{2\mu-1}\hspace{.075cm} _{2}F_{3}\left(1,-\frac{1}{2\mu}+1;-\frac{1}{2\mu}+2,\frac{3}{2},2;-\frac{\lambda^2 }{4 x^{2\mu}}\right)+C, \mu>1.
\label{eq51}
\end{align}
Making the substitution $u=1/x$ and applying Theorem \ref{teoen:5}, then for $\mu=1$, we have 
\begin{align}
\int\cos{\left(\frac{\lambda} {x}\right)}dx&=-\int\frac{\cos{\left(\lambda u\right)}}{u^2}du
=\frac{1}{2\lambda u^{2}}+\frac{\lambda}{2}\ln|u|-\frac{\lambda}{6}\left(\frac{\lambda u}{4}\right)^2\hspace{.075cm} _{2}F_{3}\left(1,1;2,\frac{5}{2},3;-\frac{\lambda^2 x^{2}}{4}\right)
\nonumber\\ &=\frac{x^2}{2\lambda}-\frac{\lambda}{2}\ln|x|-\frac{\lambda}{6}\left(\frac{\lambda}{4x}\right)^2\hspace{.075cm} _{2}F_{3}\left(1,1;2,\frac{5}{2},3;-\frac{\lambda^2}{4  x^{2}}\right)+C.
\label{eq51-1}
\end{align}
\end{enumerate}
%\label{ex:4}
%\end{exaen}
%\\[3ex]

\section{Evaluation of exponential (Ei) and logarithmic (Li) integrals}\label{sec:5}
\begin{teoen}
If $\beta\ge1$, then for any constant $\lambda$,
\begin{equation}
\int \frac{e^{\lambda x^\beta}}{x}dx=\ln{|x|}+\frac{\lambda x^\beta}{\beta}\hspace{.075cm}   _{2}F_{2}(1,1;2,2;\lambda x^\beta)+C,
\label{eq52}
\end{equation}
and
\begin{equation}
\lambda x^\beta\hspace{.075cm}   _{2}F_{2}(1,1;2,2;\lambda x^\beta)\sim-2+\frac{e^{\lambda x^\beta}}{\lambda x^\beta}, \hspace{.12cm} |x|\gg1.
\label{eq53}
\end{equation}
\label{teoen:6}
\end{teoen}

\proofen
\begin{align}
\int \frac{e^{\lambda x^\beta}}{x}dx&=\int\frac{1}{x}\sum\limits_{n=0}^{\infty}\frac{(\lambda x^\beta)^n}{n!}dx=\int\frac{dx}{x}+\int\sum\limits_{n=1}^{\infty}\frac{\lambda^n x^{\beta n-1}}{n!}dx
=\ln{|x|}+\sum\limits_{n=1}^{\infty}\frac{\lambda^n}{n!}\int x^{\beta n-1} dx
\nonumber\\&=\ln{|x|}+\sum\limits_{n=1}^{\infty}\frac{\lambda^n}{n!}\frac{x^{\beta n}}{\beta n}
=\ln{|x|}+\sum\limits_{n=0}^{\infty}\frac{\lambda^{n+1}}{(n+1)!}\frac{x^{\beta n+\beta}}{\beta n+\beta}
\nonumber\\&=\ln{|x|}+\frac{\lambda x^\beta}{\beta}\sum\limits_{n=0}^{\infty}\frac{\Gamma(n+1)}{\Gamma(n+2)\Gamma(n+2)}(\lambda x^\beta)^{n}+C
\nonumber\\&=\ln{|x|}+\frac{\lambda x^\beta}{\beta}\sum\limits_{n=0}^{\infty}\frac{(1)_n (1)_n}{(2)_n (2)_n}\frac{(\lambda x^\beta)^{n}}{n!}+C
=\ln{|x|}+\frac{\lambda x^\beta}{\beta}\hspace{.075cm}   _{2}F_{2}(1,1;2,2;\lambda x^\beta)+C.
\label{eq54}
\end{align}

To derive the asymptotic expression of $\lambda x^\beta \hspace{.075cm} _{2}F_{2}(1,1;2,2;\lambda x^\beta)$, $|x|\gg 1$, we use the asymptotic expression of the hypergeometric function ${}_2F_2\left(a_1,a_2;b_1,b_2;z\right)$ for $|z|\gg 1$, where $z\in \mathbb{C}$, and $a_1,a_2,b_1$ and $b_2$ are constants. It can be obtained using formulas 16.11.1, 16.11.2 and 16.11.7 in \cite{N} and is given by
\begin{align}
&{}_2F_2\left(a_1,a_2;b_1,b_2;z\right)\nonumber\\&=\frac{\Gamma(b_1)\Gamma(b_2)}{\Gamma(a_2)}(ze^{\pm i\pi})^{-a_1}\left\{\sum\limits_{n=0}^{R-1}\frac{(a_1)_n \Gamma(a_1-a_2-n)}{\Gamma(b_1-a_1-n)\Gamma(b_2-a_1-n)_n}\frac{(ze^{\pm i\pi})^{-n}}{n!}+O(|z|^{-R})\right\}\nonumber\\&+
\frac{\Gamma(b_1)\Gamma(b_2)}{\Gamma(a_1)}(ze^{\pm i\pi})^{-a_2}\left\{\sum\limits_{n=0}^{R-1}\frac{(a_2)_n \Gamma(a_2-a_1-n)}{\Gamma(b_1-a_2-n)\Gamma(b_2-a_2-n)_n}\frac{(ze^{\pm i\pi})^{-n}}{n!}+O(|z|^{-R})\right\}
\nonumber\\&
+\frac{\Gamma(b_1)\Gamma(b_2)}{\Gamma(a_1)\Gamma(a_2)} e^{z}z^{a_1+a_2-b_1-b_2}\left\{\sum\limits_{n=0}^{S-1}\frac{\mu_n}{2^n}z^{-n}+O(|z|^{-S})\right\},
\label{eq55}
\end{align}
where the coefficient $\mu_n$ is given by formula 16.11.4. And the upper or lower signs are chosen according as $z$ lies in the upper (above the real axis) or lower half-plane (below the real axis).

Setting $z=\lambda x^\beta, a_1=1, a_2=1, b_1=2$ and $b_2=2$ in (\ref{eq55}) yields
\begin{equation}
_{2}F_{2}(1,1;2,2;\lambda x^\beta)\sim\frac{-2}{\lambda x^\beta}+\frac{e^{\lambda x^\beta}}{\lambda^2 x^{2\beta}}, \hspace{.12cm} |x|\gg1.
\label{eq56}
\end{equation}
Hence,
\begin{equation}
\lambda x^\beta \hspace{.075cm} _{2}F_{2}(1,1;2,2;\lambda x)\sim-2+\frac{e^{\lambda x^\beta}}{\lambda x^\beta}, \hspace{.12cm} |x|\gg1.
\label{eq57}
\end{equation}
This ends the proof.
\hfill $\square$

\emph{Example~5.}
The random attenuation capacity of a channel or fading capacity \cite{MM} can now be evaluated in terms of the natural logarithm $\ln $ and the hypergeometric function $_{2}F_{2}$ as 
\begin{equation}
C_{\mbox{fading}}=E[\log_2(1+P|H|^2)]=\frac{1}{\ln{2}}e^{\frac{1}{p}}
\left[E_{1,1}\left(\infty\right)-E_{1,1}\left(\frac{1}{p}\right)\right]=\frac{1}{\ln{2}}e^{\frac{1}{p}}\left[\ln{P}+\frac{1}{P}\ _{2}F_{2}\left(1,1;2,2;-\frac{1}{P}\right)\right]
\label{application1-1}
\end{equation}
\\[5ex]
\emph{Example~6.}
%\begin{exaen}
One can now evaluate $\int e^{\lambda e^{\beta x}}dx$ in terms of $_{2}F_{2}$ using the substitution $u=e^x$, and obtain
\begin{align}
\int e^{\lambda e^{\beta x} }dx=\int \frac{e^{\lambda u^\beta}}{u} du&=\ln{u}+ \frac{\lambda u^\beta}{\beta}\hspace{.075cm}   _{2}F_{2}(1,1;2,2; \lambda u^\beta)+C=x+ \frac{\lambda e^{\beta x}}{\beta}\hspace{.075cm}   _{2}F_{2}(1,1;2,2; \lambda e^{\beta x})+C.
\label{eq58}
\end{align}
%\label{ex:5}
%\end{exaen}
\\[6ex]

\begin{teoen} The logarithmic integral is given by
\begin{equation}
\text{Li}=\int\limits_\mu^x\frac{dt}{\ln{t}}=\ln{\left(\frac{\ln{x}}{\ln{\mu}}\right)}+\ln{x}\hspace{.075cm} _{2}F_{2}(1,1;2,2;\ln{x})-
\ln{\mu}\hspace{.075cm} _{2}F_{2}(1,1;2,2;\ln{\mu}), \mu>1.
\label{eq59}
\end{equation}
And for $x\gg\mu$,
\begin{equation}
\text{Li}=\int\limits_{\mu}^{x}\frac{dt}{\ln{t}}\sim \frac{x}{\ln{x}}+\ln{\left(\frac{\ln{x}}{\ln{\mu}}\right)}-2-
\ln{\mu}\hspace{.075cm} _{2}F_{2}(1,1;2,2;\ln{\mu}).
\label{eq60}
\end{equation}
\label{teoen:7}
\end{teoen}

\proofen
Making the substitution $u=\ln{x}$ and using (\ref{eq50}) gives
\begin{align}
\int\limits_{\mu}^{x}\frac{dx}{\ln{x}}&=\int\limits_{\ln\mu}^{\ln{x}}\frac{e^u}{u}du=\left[\ln{u}+u\hspace{.075cm} _{2}F_{2}(1,1;2,2;u)\right]_{\ln\mu}^{\ln{x}}
\nonumber\\ & =\ln{\left(\frac{\ln{x}}{\ln\mu}\right)}+\ln{x}\hspace{.075cm} _{2}F_{2}(1,1;2,2;\ln{x})-
\ln{\mu}\hspace{.075cm} _{2}F_{2}(1,1;2,2;\ln{\mu}).
\label{eq61}
\end{align}
Now setting $z=\ln{x}, a_1=1$, $a_2=1$, $b_1=2$ and $b_2=2$ in (\ref{eq53}) or in (\ref{eq51}) yields
\begin{equation}
_{2}F_{2}(1,1;2,2;\ln{x})\sim\frac{-2}{\ln{x}}+\frac{x}{(\ln{x})^2}, \hspace{.2cm} x\gg1.
\label{eq62}
\end{equation}
This gives
\begin{equation}
\ln{x}\hspace{.075cm} _{2}F_{2}(1,1;2,2;\ln{x})\sim-2+\frac{x}{\ln{x}}, \hspace{.2cm} x\gg1.
\label{eq63}
\end{equation}
Hence for $x\gg\mu$,
\begin{equation}
\text{Li}=\int\limits_{\mu}^{x}\frac{dt}{\ln{t}}\sim \frac{x}{\ln{x}}+\ln{\left(\frac{\ln{x}}{\ln{\mu}}\right)}-2-
\ln{\mu}\hspace{.075cm} _{2}F_{2}(1,1;2,2;\ln{\mu}).
\label{eq64}
\end{equation}
\hfill $\square$

We importantly note that Theorem \ref{teoen:7} adds the term $\ln{\left(\frac{\ln{x}}{\ln{\mu}}\right)}-2-
\ln{\mu}\hspace{.075cm} _{2}F_{2}(1,1;2,2;\ln{\mu})$ to the known asymptotic expression of the logarithmic integral in mathematical litterature, $\text{Li}\sim {x}/{\ln{x}}$ \cite{AS,N}. And this term is negligible if $x\sim O(10^6)$ or higher.

We can now slightly improve the prime number Theorem \cite{H} as following,
\begin{coren} Let $\pi(x)$ denotes the number of primes small than or equal to $x$ and $\mu>1$. Then for $x\gg\mu$,
\begin{equation}
\pi(x)-\frac{x}{\ln{x}}\sim\ln{\left(\frac{\ln{x}}{\ln{\mu}}\right)}-2-
\ln{\mu}\hspace{.075cm} _{2}F_{2}(1,1;2,2;\ln{\mu}).
\label{eq60-1}
\end{equation}
\label{coren:2}
\end{coren}
The proof follows directly from equation (\ref{eq60}) in Theorem \ref{teoen:7}.
 
\emph{Example~7.}
One can now evaluate $\int \ln{(\ln{x})} dx$ using integration by parts.
\begin{align}
\int \ln{(\ln{x})} dx&=x\ln{(\ln{x})}-\int \frac{1}{\ln{x}} dx
\nonumber \\ &=x\ln{(\ln{x})}-\ln{(\ln{x})}-\ln{x}\hspace{.075cm} _{2}F_{2}(1,1;2,2;\ln{x})+C.
\label{eq65}
\end{align}
%\label{ex:6}
%\end{exaen}
\\[7ex]
\begin{teoen}
For $\beta\ge1$ and $\alpha<\beta+1$, we have
\begin{equation}
\text{Ei}_{\beta,\alpha}=\int \frac{e^{\lambda x^\beta}}{\lambda x^\alpha}dx=\frac{1}{\lambda}\frac{x^{1-\alpha}}{1-\alpha}+\frac{x^{\beta-\alpha+1}}{\beta-\alpha+1}\hspace{.075cm}   _{2}F_{2}\left(1,-\frac{\alpha}{\beta}+\frac{1}{\beta}+1;2,-\frac{\alpha}{\beta}+\frac{1}{\beta}+2;\lambda x^\beta\right)+C,
\label{eq66}
\end{equation}
and for $|x|\gg1$,
\begin{multline}
\frac{\lambda x^{\beta-\alpha+1}}{\beta-\alpha+1} \hspace{.075cm}_{2}F_{2}\left(1,-\frac{\alpha}{\beta}+\frac{1}{\beta}+1;2,-\frac{\alpha}{\beta}+\frac{1}{\beta}+2;\lambda x\right)\\
\sim \frac{\lambda}{\beta} \Gamma\left(-\frac{\alpha}{\beta}+\frac{1}{\beta}+1\right)\left(-\frac{1}{\lambda }\right)^{-\frac{\alpha}{\beta}+\frac{1}{\beta}+1}-\frac{x^{-\alpha+1}}{\beta}
+\frac{1}{\lambda\beta}\frac{e^{\lambda x^\beta}}{x^{\beta+\alpha-1}}.
\label{eq67}
\end{multline}
We also have,
\begin{equation}
\text{Ei}_{\beta,\beta+1}=\int \frac{e^{\lambda x^\beta}}{\lambda x^{\beta+1}}dx=-\frac{1}{\beta x^\beta}+\ln(|x|)+\frac{\lambda x^{\beta}}{2\beta}\hspace{.075cm}   _{2}F_{2}\left(1,1;2,2;\lambda x^\beta\right)+C.
\label{eq66-1}
\end{equation} 
\label{teoen:8}
\end{teoen}

\proofen
\begin{align}
\text{Ei}_{\beta,\alpha}&=\int \frac{e^{\lambda x^\beta}}{\lambda x^\alpha}dx=\int\frac{1}{\lambda x^\alpha}\sum\limits_{n=0}^{\infty}\frac{(\lambda x^\beta)^n}{n!}dx
=\frac{1}{\lambda}\int\frac{dx}{x^\alpha}+=\frac{1}{\lambda x^\alpha}\int\sum\limits_{n=1}^{\infty}\frac{(\lambda x^\beta)^n}{n!}dx
\nonumber\\ &=\frac{1}{\lambda}\frac{x^{1-\alpha}}{1-\alpha}+\frac{1}{\lambda}\sum\limits_{n=1}^{\infty}\frac{\lambda^n}{n!}\int x^{\beta n-\alpha}dx
=\frac{1}{\lambda}\frac{x^{1-\alpha}}{1-\alpha}+\sum\limits_{n=0}^{\infty}\frac{\lambda^{n}}{(n+1)!}\frac{ x^{\beta n+\beta-\alpha+1}}{\beta n+\beta-\alpha+1}+C
\nonumber\\ &=\frac{1}{\lambda}\frac{x^{1-\alpha}}{1-\alpha}+\frac{x^{\beta-\alpha+1}}{\beta}\sum\limits_{n=0}^{\infty}\frac{\Gamma\left(n-\frac{\alpha}{\beta}+\frac{1}{\beta}+1\right)}{\Gamma(n+2)\Gamma\left(n-\frac{\alpha}{\beta}+\frac{1}{\beta}+2\right)}
\left(\lambda x^{\beta}\right)^n+C
\nonumber\\ &=\frac{1}{\lambda}\frac{x^{1-\alpha}}{1-\alpha}+\frac{x^{\beta-\alpha+1}}{\beta-\alpha+1}
\sum\limits_{n=0}^{\infty}\frac{(1)n\left(-\frac{\alpha}{\beta}+\frac{1}{\beta}+1\right)_n}{(2)_n\left(-\frac{\alpha}{\beta}+\frac{1}{\beta}+2\right)_n}\frac{ \left( x^{\beta}\right)^n}{n!}+C
\nonumber\\ &=\frac{1}{\lambda}\frac{x^{1-\alpha}}{1-\alpha}+\frac{ x^{\beta-\alpha+1}}{\beta-\alpha+1}\hspace{.075cm}   _{2}F_{2}\left(1,-\frac{\alpha}{\beta}+\frac{1}{\beta}+1;2,-\frac{\alpha}{\beta}+\frac{1}{\beta}+2;\lambda x^\beta\right)+C.
\label{eq68}
\end{align}
Now setting $a_1=1, a_2=-\frac{\alpha}{\beta}+\frac{1}{\beta}+1, b_1=2, b_2=-\frac{\alpha}{\beta}+\frac{1}{\beta}+2$ and $z=\lambda x^\beta$ in (\ref{eq55}) gives,
\begin{multline}
_{2}F_{2}\left(1,-\frac{\alpha}{\beta}+\frac{1}{\beta}+1;2,-\frac{\alpha}{\beta}+\frac{1}{\beta}+2;\lambda x^\beta\right)\\
\sim -\left(-\frac{\alpha}{\beta}+\frac{1}{\beta}+1\right)\frac{1}{\lambda x^\beta}
+\Gamma\left(-\frac{\alpha}{\beta}+\frac{1}{\beta}+2\right)\left(\frac{1}{\lambda x^\beta}\right)^{-\frac{\alpha}{\beta}+\frac{1}{\beta}+1}+\frac{e^{\lambda x^\beta}}{\lambda^2 x^{2\beta}}.
\label{eq69}
\end{multline}
Hence, multiplying (\ref{eq69}) with $\frac{\lambda x^{\beta-\alpha+1}}{\beta-\alpha+1}$ gives (\ref{eq67}). The proof of  (\ref{eq66-1}) is similar to that of (\ref{eq45}).
\hfill $\square$

\begin{teoen}
For any constants $\alpha$, $\beta$ and $\lambda$,
\begin{align}
&_{1}F_{2}\left(-\frac{\alpha}{2\beta}+\frac{1}{2\beta}+\frac{1}{2};-\frac{\alpha}{2\beta}+\frac{1}{2\beta}+\frac{3}{2},\frac{3}{2};-\frac{\lambda^2 x^{2\beta}}{4}\right)=\nonumber\\ & \frac{1}{2}\Bigl[{}_{2}F_{2}\left(1,-\frac{\alpha}{\beta}+\frac{1}{\beta}+1;2,-\frac{\alpha}{\beta}+\frac{1}{\beta}+2;i\lambda x^\beta\right)+{}_{2}F_{2}\left(1,-\frac{\alpha}{\beta}+\frac{1}{\beta}+1;2,-\frac{\alpha}{\beta}+\frac{1}{\beta}+2;-i\lambda x^\beta\right)\Bigr],
\label{eq70}
\end{align}
%\item
or
\begin{align}
&_{1}F_{2}\left(-\frac{\alpha}{2\beta}+\frac{1}{2\beta}+\frac{1}{2};-\frac{\alpha}{2\beta}+\frac{1}{2\beta}+\frac{3}{2},\frac{3}{2};\frac{\lambda^2 x^{2\beta}}{4}\right)=\nonumber\\&\frac{1}{2}\Bigl[{}_{2}F_{2}\left(1,-\frac{\alpha}{\beta}+\frac{1}{\beta}+1;2,-\frac{\alpha}{\beta}+\frac{1}{\beta}+2;\lambda x^\beta\right)+{}_{2}F_{2}\left(1,-\frac{\alpha}{\beta}+\frac{1}{\beta}+1;2,-\frac{\alpha}{\beta}+\frac{1}{\beta}+2;-\lambda x^\beta\right)\Bigr].
\label{eq71}
\end{align}
%\end{enumerate}
\label{teoen:9}
\end{teoen}

\proofen
%\begin{enumerate}
%\item
Using Theorem \ref{teoen:8}, we obtain
\begin{multline}
\int\frac{\sin{(\lambda  x^\beta)}}{x^\alpha} dx=\frac{1}{2i}\int  \frac{e^{i\lambda x^\beta}-e^{-i\lambda x^\beta}}{x^\alpha}dx
\\=\frac{1}{2}\frac{\lambda x^{\beta-\alpha+1}}{\beta-\alpha+1}\Bigl[ \hspace{.075cm}_{2}F_{2}\left(1,-\frac{\alpha}{\beta}+\frac{1}{\beta}+1;2,-\frac{\alpha}{\beta}+\frac{1}{\beta}+2;i\lambda x^\beta\right)\\+ \hspace{.075cm}_{2}F_{2}\left(1,-\frac{\alpha}{\beta}+\frac{1}{\beta}+1;2,-\frac{\alpha}{\beta}+\frac{1}{\beta}+2;-i\lambda x^\beta\right)\Bigr]+C.
\label{eq72}
\end{multline}
Hence, comparing (\ref{eq23}) with (\ref{eq72}) gives (\ref{eq70}). Or on the other hand,
\begin{multline}
2\int\frac{\sinh{(\lambda  x^\beta)}}{x^\alpha} dx=\int  \frac{e^{\lambda x^\beta}-e^{-\lambda x^\beta}}{x^\alpha}dx=\frac{\lambda x^{\beta-\alpha+1}}{\beta-\alpha+1}\times
\\\Bigl[ \hspace{.075cm}_{2}F_{2}\left(1,-\frac{\alpha}{\beta}+\frac{1}{\beta}+1;2,-\frac{\alpha}{\beta}+\frac{1}{\beta}+2;\lambda x^\beta\right)+ _{2}F_{2}\left(1,-\frac{\alpha}{\beta}+\frac{1}{\beta}+1;2,-\frac{\alpha}{\beta}+\frac{1}{\beta}+2;-\lambda x^\beta\right)\Bigr]+C.
\label{eq73}
\end{multline}
Hence, comparing (\ref{eq28}) with (\ref{eq75}) gives (\ref{eq73}).

%\end{enumerate}

\hfill $\square$

\begin{teoen}
For any constants $\alpha$, $\beta$ and $\lambda$,
\begin{multline}
\frac{i x^{2\beta-\alpha+1}}{2\beta-\alpha+1}\hspace{.075cm} _{2}F_{3}\left(1,-\frac{\alpha}{2\beta}+\frac{1}{2\beta}+1;-\frac{\alpha}{2\beta}+\frac{1}{2\beta}+2,\frac{3}{2},2;-\frac{\lambda^2 x^{2\beta}}{4}\right)=
\frac{x^{\beta-\alpha+1}}{\beta-\alpha+1}\times\\\Bigl[\hspace{.075cm}   _{2}F_{2}\left(1,-\frac{\alpha}{\beta}-\frac{1}{\beta}+1;2,-\frac{\alpha}{\beta}+\frac{1}{\beta}+2;i\lambda x^\beta\right)
 -\hspace{.075cm}_{2}F_{2}\left(1,-\frac{\alpha}{\beta}+\frac{1}{\beta}+1;2,-\frac{\alpha}{\beta}+\frac{1}{\beta}+2;-i\lambda x^\beta\right)\Bigr].
\label{eq74}
\end{multline}
Or, \begin{multline}
\frac{x^{2\beta-\alpha+1}}{2\beta-\alpha+1}\hspace{.075cm} _{2}F_{3}\left(1,-\frac{\alpha}{2\beta}+\frac{1}{2\beta}+1;-\frac{\alpha}{2\beta}+\frac{1}{2\beta}+2,\frac{3}{2},2;\frac{\lambda^2 x^{2\beta}}{4}\right)\\=
\frac{x^{\beta-\alpha+1}}{\beta-\alpha+1}\Bigl[\hspace{.075cm}   _{2}F_{2}\left(1,-\frac{\alpha}{\beta}-\frac{1}{\beta}+1;2,-\frac{\alpha}{\beta}+\frac{1}{\beta}+2;\lambda x^\beta\right)
\\  +\hspace{.075cm}_{2}F_{2}\left(1,-\frac{\alpha}{\beta}+\frac{1}{\beta}+1;2,-\frac{\alpha}{\beta}+\frac{1}{\beta}+2;-\lambda x^\beta\right)\Bigr].
\label{eq75}
\end{multline}
%\end{enumerate}
\label{teoen:10}
\end{teoen}
We prove Theorem \ref{teoen:10} as Theorem  \ref{teoen:9} using Theorems \ref{teoen:4} and \ref{teoen:8}.

%\proofen
%Using Theorem \ref{teoen:8}, we obtain
%\begin{multline}
%2\int \frac{\cos(\lambda x^\beta)}{x^\alpha}dx=\int\frac{e^{i\lambda x^\beta}+e^{-i\lambda x^\beta}}{x^\alpha}dx
%=\frac{x^{1-\alpha}}{1-\alpha}\\+\frac{i\lambda x^{\beta-\alpha+1}}{\beta-\alpha+1} \Bigl[ \hspace{.075cm}_{2}F_{2}\left(1,-\frac{\alpha}{\beta}+\frac{1}{\beta}+1;2,-\frac{\alpha}{\beta}+\frac{1}{\beta}+2;i\lambda x^\beta\right)\\- \hspace{.075cm}_{2}F_{2}\left(1,-\frac{\alpha}{\beta}+\frac{1}{\beta}+1;2,-\frac{\alpha}{\beta}+\frac{1}{\beta}+2;-i\lambda x^\beta\right)\Bigr]+C.
%\label{eq76}
%\end{multline}
%Hence, comparing (\ref{eq39}) with (\ref{eq76}) gives (\ref{eq74}). Or on the other hand,
%\begin{multline}
%2\int \frac{\cosh(\lambda x^\beta)}{x^\alpha}dx=\int\frac{e^{\lambda x^\beta}+e^{-\lambda x^\beta}}{x^\alpha}dx
%=\frac{x^{1-\alpha}}{1-\alpha}\\+\frac{\lambda x^{\beta-\alpha+1}}{\beta-\alpha+1} \Bigl[ \hspace{.075cm}_{2}F_{2}\left(1,-\frac{\alpha}{\beta}+\frac{1}{\beta}+1;2,-\frac{\alpha}{\beta}+\frac{1}{\beta}+2;\lambda x^\beta\right)\\- \hspace{.075cm}_{2}F_{2}\left(1,-\frac{\alpha}{\beta}+\frac{1}{\beta}+1;2,-\frac{\alpha}{\beta}+\frac{1}{\beta}+2;-\lambda x^\beta\right)\Bigr]+C.
%\label{eq77}
%\end{multline}
%Hence, comparing (\ref{eq44}) with (\ref{eq77}) gives (\ref{eq75}).
%\hfill $\square$

\section{Conclusion}\label{sec:6}
\setcounter{equation}{0}
$\text{Si}_{\beta,\alpha}=\int [\sin{(\lambda x^\beta)}/{(\lambda x^\alpha)}] dx, \beta\ge1, \alpha\le\beta+1$, and $\text{Ci}_{\beta,\alpha}=\int [\cos{(\lambda x^\beta)}/{(\lambda x^\alpha)}] dx,\beta\ge1, \alpha\le 2\beta+1$, were expressed in terms of the hypergeometric functions ${}_1F_2$ and ${}_2F_3$ respectively, and their asymptotic expressions for $|x|\gg1$ were obtained (see Theorems \ref{teoen:1},\ref{teoen:2}, \ref{teoen:3}, \ref{teoen:4} and \ref{teoen:5}). Once derived, formulas for the hyperbolic sine and hyperbolic cosine integrals were readily deduced from those of the sine and cosine integrals.

On the other hand, the exponential integral $\text{Ei}_{\beta,\alpha}=\int (e^{\lambda x^\beta}/x^\alpha)dx, \beta\ge1, \alpha\le \beta+1$,  and the logarithmic integral $\int dx/\ln{x}$ were expressed in terms of the hypergeometric function ${}_2F_2$, and their asymptotic expressions for $|x|\gg1$ were also obtained (see Theorems \ref{teoen:6}, \ref{teoen:7} and  \ref{teoen:8}). Therefore, their corresponding definite integrals can now be evaluated using the FTC rather than using numerical integration.

Using the Euler and hyperbolic identities $\text{Si}_{\beta,\alpha}$ and $\text{Ci}_{\beta,\alpha}$ were expressed in terms of $\text{Ei}_{\beta,\alpha}$. And hence, some expressions of the hypergeometric functions ${}_1F_2$ and ${}_2F_3$ in terms of ${}_2F_2$ were derived (see Theorems \ref{teoen:9} and  \ref{teoen:10}).

The evaluation of the logarithmic integral $\int dx/\ln{x}$ in terms of the function ${}_2F_2$ and its asymptotic expression ${}_2F_2$ for $|x|\gg1$ allowed us to add the term $\ln{\left(\frac{\ln{x}}{\ln{\mu}}\right)}-2-\ln{\mu}\hspace{.075cm} _{2}F_{2}(1,1;2,2;\ln{\mu}), \mu>1,$ to the known asymptotic expression of the logarithmic integral, which is $\text{Li}=\int_{2}^{x} dt/\ln{t}\sim {x}/{\ln{x}}$ \cite{AS,N}, so that it is given by $\text{Li}=\int_{\mu}^{x}{dt}/{\ln{t}}\sim {x}/{\ln{x}}+\ln{\left(\frac{\ln{x}}{\ln{\mu}}\right)}-2-\ln{\mu}\hspace{.075cm} _{2}F_{2}(1,1;2,2;\ln{\mu})$ in Theorem \ref{teoen:7}. Beside, this leads to Corollary \ref{coren:2} which is an improvement of the prime number Theorem \cite{H}.

In addition, other non-elementary integrals  which can be written in terms of $\text{Ei}_{\beta,1}$ and  $\int dx/\ln{x}$ and then evaluated were given as examples. For instance, using substitution, the $\int e^{\lambda e^ {\beta x}}dx$ was written in terms of $\text{Ei}_{\beta,1}$ and therefore evaluated in terms of $_{2}F_{2}$, and using integration by parts, the non-elementary integral $\int \ln(\ln{x}) dx$ was written in terms of $\int dx/\ln{x}$ and therefore evaluated in terms of $_{2}F_{2}$.

%The goal of this paper was to evaluate the above integrals in terms of elementary functions and special functions for which properties are well known so that the FTC can be used in evaluating their definite integrals. So the cases $\text{Si}_{\beta,\alpha}=\int [\sin{(\lambda x^\beta)}/{(\lambda x^\alpha)}] dx, \beta\ge1, \alpha>\beta+1$, and $\text{Ci}_{\beta,\alpha}=\int [\cos{(\lambda x^\beta)}/{(\lambda x^\alpha)}] dx,\beta\ge1, \alpha> 2\beta+1$ and $\text{Ei}_{\beta,\alpha}=\int (e^{\lambda x^\beta}/x^\alpha)dx, \beta\ge1, \alpha> \beta+1$, that may involve series whose properties are not known yet and requiring the use of numerical methods, as already mentioned, will be treated in Part 2.  

\begin{Biblioen}

\bibitem{AS}{\bf Abramowitz~M., Stegun~I.A.} Handbook of mathematical functions with formulas,
graphs and mathematical tables. National Bureau of Standards,1964. 1046~p.
\bibitem{B}{\bf Billingsley~P.} Probability and measure. Wiley series in Probability and Mathematical
Statistics, 3rd Edition, 1995. 608~p.
\bibitem{H}{\bf Hoffman~P.}The man who loved only numbers. New York Hyperion Books, 1998. 227~p.
\bibitem{K}{\bf Krantz~S.G.} Handbook of Complex variables. Boston: MA Birk¨ausser, 1999. 290~p.
\bibitem{L}{\bf Lebedev~ N.N.} Special functions and their applications. Prentice-Hall Inc., Englewood Cliffs, N.J., 1965.
\bibitem{MZ}{\bf Marchisotto~E.A., Zakeri~G.-A.} An invitation to integration in finite terms// College
Math. J., 1994. Vol.~25, no~4. P.~295--308. DOI:~10.2307/2687614
\bibitem{NV}{\bf Nijimbere~V.} Evaluation of the non-elementary integral $\int e^{\lambda x^\alpha} dx$, $\alpha\ge2$, and other related integrals// Ural Math. J., 2017. Vol~3, no.~2. P.~130--142. DOI:~10.15826/umj.2017.2.014
\bibitem{NV2}{\bf Nijimbere~V.} Evaluation of some non-elementary integrals involving sine, cosine, exponential and logarithmic integrals: Part II// Ural Math. J., 2017. Accepted for publication.
\bibitem{N} NIST Digital Library of Mathematical Functions. \url{http://dlmf.nist.gov/}
\bibitem{R}{\bf Rosenlicht~M.} Integration in finite terms// Amer. Math. Monthly, 1972. Vol~79, no.~9. P.~963--972. 
DOI:~10.2307/2318066
\bibitem{MM}{\bf Simon~M.K., Alouini~M.‐S.} Digital Communication over Fading Channels. John Wiley \& Sons, Inc., 2005. 900~p.
\end{Biblioen}

 \end{document}